# Two-Dimensional Cruise Control of Autonomous Vehicles on Lane-Free Roads

**Iasson Karafyllis[*], Dionysis Theodosis[**] and Markos Papageorgiou[**]**

[*]Dept. of Mathematics, National Technical University of Athens,
Zografou Campus, 15780, Athens, Greece,
emails: iasonkar@central.ntua.gr , iasonkaraf@gmail.com

[**] Dynamic Systems and Simulation Laboratory,
Technical University of Crete, Chania, 73100, Greece
(emails: dtheodosis@dssl.tuc.gr , markos@dssl.tuc.gr)

**Abstract**

In this paper, we design decentralized control strategies for the two-dimensional movement of autonomous vehicles on lane-free roads. The bicycle kinematic model is used to model the dynamics of the vehicles, and each vehicle determines its control input based only on its own speed and on the distance from other (adjacent) vehicles and the boundary of the road. Potential functions and Barbălat's lemma are employed to prove the following properties, which are ensured by the proposed controller: (i) the vehicles do not collide with each other or with the boundary of the road; (ii) the speeds of all vehicles are always positive, i.e., no vehicle moves backwards at any time; (iii) the speed of all vehicles remain below a given speed limit; (iv) all vehicle speeds converge to a given longitudinal speed set-point; and (v) the accelerations, lateral speeds, and orientations of all vehicles tend to zero. The efficiency of the proposed 2-D cruise controllers is illustrated by means of numerical examples.

**Keywords:** Cruise Control, Lane-Free Roads, Autonomous Vehicles, Lyapunov functions.

## 1. Introduction

Vehicle automation has made tremendous advances in the last decades, and the path to full automation of vehicles in a foreseeable future seems more than likely. An initial stage of vehicle automation is the standard cruise control system which maintains the speed of the vehicle at a desired value to assist the driver. These systems have meanwhile evolved to Adaptive Cruise Control (ACC) systems, which automatically adjust the speed to maintain certain distance from a front vehicle or to maintain a desired speed. Recent advances of communication technologies have also been used in vehicle automation to develop Cooperative ACC systems (CACC) so that vehicles are wirelessly connected and can therefore communicate to maintain shorter inter-vehicle distances, thus increasing capacity of the roads. It has been shown that communication between vehicles has the potential of increasing their safety, reduce congestion and traffic accidents and improve traffic flow on highways ([4], [15], [24]). Both ACC and CACC systems have been extensively studied in the literature (see for instance [9], [12], [15], [21], [25], [28]).

The vast majority of research effort is focused on studying lane-based traffic models, where vehicles abide to a lane discipline, which increases traffic safety, as it simplifies the task of manual driving. Indeed, all control strategies for ACC and CACC systems are developed based on information from the vehicle directly in front or behind (see for instance [9], [12], [15], [21], [25] and references therein). Apart from the car-following task, another necessary driving task is lane-changing, which is a more complex and riskier maneuver, since the driver needs to look for an available gap on the target lane and estimate the speeds of many adjacent vehicles quasi-simultaneously. Modeling lane changes and two-dimensional movement on lane-based roads is a complicated problem, and various approaches have been considered, see for instance [6], [22], [30]. Lane changing is known to affect both capacity and safety [22], [30].

Recently, launched by [19], new principles and research directions were proposed for autonomous vehicles operating on lane-free roads ([3], [16], [19]) that may improve traffic flow and increase capacity of highways. The vehicles move on the two-dimensional surface of the road without obeying to a lane discipline as in conventional traffic. Since connected and automated vehicles use sensors and can communicate their presence and state to other vehicles, they are suitable and more efficient in a lane-free environment where they can use their capabilities to their full extent. For the lane-free concept, only a few models have been proposed that can describe vehicle movement on lane-free roads, driven by human drivers; see [1, 10, 17]. These approaches are based on linear systems theory and traditional longitudinal car-following models, which, however, do not guarantee: (i) collision avoidance with other vehicles or the boundary of the road, (ii) positivity of speeds, and (iii) speeds within road speed limits. In addition to the lane-free traffic, another concept that can increase the flow of vehicles on a road is the associated concept of 'nudging' (see [19]). Nudging implies a virtual force that vehicles apply to the vehicles in front of them, and it has been shown that nudging can increase the flow in a ring-road and can have a strong stabilizing effect; see [11] and references therein.

In this paper, we consider identical autonomous vehicles described by the bicycle kinematic model. The bicycle kinematic model is selected because it is able to capture the non-holonomic constraints of the actual vehicle (see [20], [21]). We design a family of novel decentralized controllers for the safe operation of the vehicles on lane-free roads. The main features of the proposed approach are:

- The proposed controllers are fully decentralized, and each vehicle only has access to the distance from the boundaries of the road and the distance from adjacent vehicles and does not require any information or estimates of relative speeds or relative orientation;
- the vehicles do not collide with each other or with the boundary of the road;
- the speeds of all vehicles are always positive and remain below a given speed limit;
- all vehicle speeds converge to a given longitudinal speed set-point; and
- the accelerations, lateral speeds, and orientations of all vehicles tend to zero.

To our knowledge, a model that describes vehicle movement in two-dimensional lane-free road and captures all these properties simultaneously is not available in the literature. To better exploit the lane-free two-dimensional road surface, we adopt elliptic metrics which approximate efficiently the physical dimensions of a vehicle and allow to explicitly determine the minimum safety distance between vehicles to avoid collisions. By appropriately selecting the "elliptical" distances, we can regulate both the safety distance and the number of vehicles that can be placed laterally on the road. To avoid collisions between vehicles and with the boundary of the road, we employ potential functions, which have been extensively used to address numerous problems in multi-agent systems such as collision avoidance, flocking, dispersion, and formation, see [5], [7], [13], [17], [22], [28]. Potential functions have also been used in certain lane-changing and lane-keeping problems in traffic control, [8], [26]. Finally, we combine Lyapunov functions with barrier functions (see [2], [12], [26]) to restrain the movement of the vehicles and exploit Barbălat's lemma [13] to address the objectives of speeds, acceleration, and orientation convergence as stated above. The main challenges that arise stem from the fact that the control systems studied in the paper evolve on specific open sets, and, in addition, various objectives and constraints must be satisfied simultaneously (positive speeds within road speed limits that converge to a specific speed set-point; vehicle movement within road boundaries; bounded orientation; and accelerations, lateral speeds, and orientations that all converge to zero).

The structure of the paper is as follows. Section 2 is devoted to the presentation of the problem formulation and the objectives of the paper. Section 3 contains the main results. Section 4 presents numerical examples to demonstrate the efficiency of the proposed decentralized cruise controllers. All proofs of the main results are provided in Section 5. Finally, some concluding remarks are given in Section 6.

**Notation.** Throughout this paper, we adopt the following notation.

* $\Re_+ := [0, +\infty)$ denotes the set of non-negative real numbers.
* By $|x|$ we denote both the Euclidean norm of a vector $x \in \Re^n$ and the absolute value of a scalar $x \in \Re$.
* Let $A \subseteq \Re^n$ be an open set. By $C^0(A, \Omega)$, we denote the class of continuous functions on $A \subseteq \Re^n$, which take values in $\Omega \subseteq \Re^m$. By $C^k(A; \Omega)$, where $k \geq 1$ is an integer, we denote the class of functions on $A \subseteq \Re^n$ with continuous derivatives of order $k$, which take values in $\Omega \subseteq \Re^m$. When $\Omega = \Re$ the we write $C^0(A)$ or $C^k(A)$. For a function $V \in C^1(A; \Re)$, the gradient of $V$ at $x \in A \subseteq \Re^n$, denoted by $\nabla V(x)$, is the row vector $\left[ \frac{\partial V}{\partial x_1}(x) \cdots \frac{\partial V}{\partial x_n}(x) \right]$. By $\nabla^2 V(x)$ we denote the Hessian matrix at $x \in A \subseteq \Re^n$ of a function $V \in C^2(A; \Re)$.
* Let $C \subset \Re^n$ be a convex set. A function $f : C \to \Re$ is convex if $f(\mu x_1 + (1-\mu) x_2) \leq \mu f(x_1) + (1-\mu) f(x_2)$ for all $x_1, x_2 \in C$ and $\mu \in [0,1]$.

## 2. Problem Description

Consider $n$ identical vehicles moving on a lane-free road of width $2a > 0$, see Figure 1. The movement of the vehicles is described by the following set of ODEs:

$$\begin{aligned}
\dot{x}_i &= v_i \cos(\theta_i) \\
\dot{y}_i &= v_i \sin(\theta_i) \\
\dot{\theta}_i &= \sigma^{-1} v_i \tan(\delta_i) \\
\dot{v}_i &= F_i
\end{aligned} \tag{2.1}$$

for $i = 1,...,n$, where $\sigma > 0$ is the length of each vehicle (a constant). Here, $(x_i, y_i)$ is the reference point of the $i$-th vehicle with $i \in \{1,...,n\}$ and is placed at the midpoint of the rear axle, with $x_i \in \Re$ being the longitudinal position and $y_i \in (-a, a)$ being the lateral position of the vehicle; $v_i$ is the speed of the $i$-th vehicle at the point $(x_i, y_i)$, $\theta_i \in \left( -\frac{\pi}{2}, \frac{\pi}{2} \right)$ is the angular orientation of the $i$-th vehicle, $\delta_i$ is the steering angle of the front wheels relative to the orientation $\theta_i$ of the $i$-th vehicle, and $F_i$ is the acceleration of the $i$-th vehicle. Model (2.1) is known as the bicycle kinematic model (see Figure 1) and has been widely used to represent vehicles due to its simplicity to capture vehicle motion in normal driving conditions, see ([14], [20], [21], [22]). To make the subsequent analysis less cumbersome, we define

$$u_i = \sigma^{-1} v_i \tan(\delta_i), \ i = 1,...,n \tag{2.2}$$

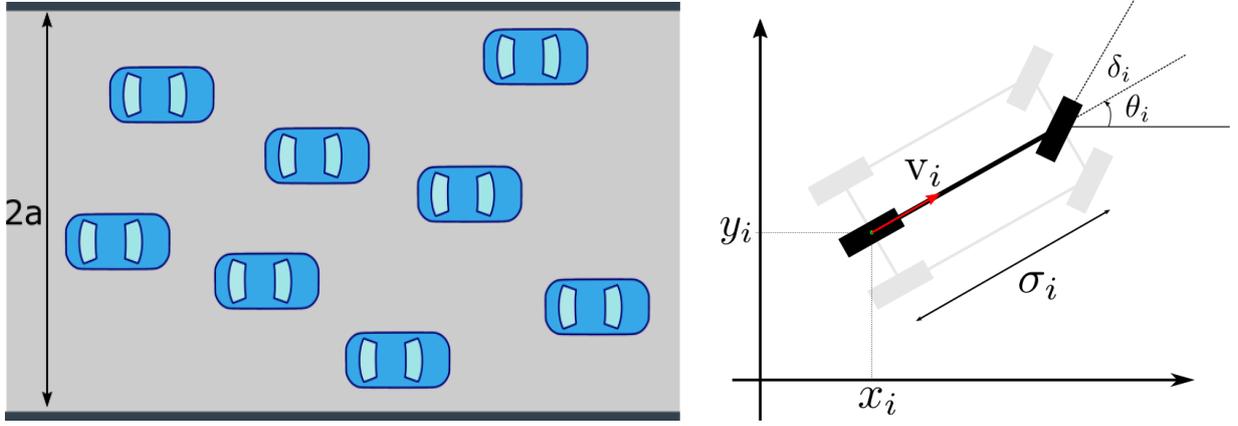

**Figure 1:** Lane-free road of width $2a > 0$ (left). Each vehicle is modeled by the bicycle kinematic model (right).

Then, model (2.1) can be written in the form

$$\begin{aligned} \dot{x}_i &= v_i \cos(\theta_i) \\ \dot{y}_i &= v_i \sin(\theta_i) \\ \dot{\theta}_i &= u_i \\ \dot{v}_i &= F_i \end{aligned} \tag{2.3}$$

for $i = 1,...,n$, where $u_i$ and $F_i$ are the control inputs. Then, $\delta_i$ can be obtained directly from (2.2) as a function of $u_i$. In what follows, we use the simpler model (2.3) instead of (2.1).

In what follows, we assume that there is no communication between the vehicles, and that the only available sensing among adjacent vehicles concerns the (elliptical) distance between vehicles, defined by

$$d_{i,j} := \sqrt{(x_i - x_j)^2 + p(y_i - y_j)^2}, \text{ for } i, j = 1,...,n \tag{2.4}$$

where $p > 0$ is a weighting factor. For $p = 1$ we obtain the standard Euclidean distance, while for larger values of $p > 1$, we have an "elliptical" metric which will allow to approximate more accurately the dimensions of a vehicle. More specifically, when vehicles are desired to maintain a certain distance between them, elliptical metrics allow to place more vehicles across the width of the road. The optimal selection of the constant $p \geq 1$ will be discussed in Section 3.

In what follows we use the notation

$$w = (x_1,...,x_n, y_1,...,y_n, \theta_1,...,\theta_n, v_1,....,v_n)' \in \Re^{4n} \tag{2.5}$$

for the stack vector of longitudinal and lateral positions, orientations and speeds off all $n$ vehicles. We assume that all vehicles operate on a lane-free road with speed limit $v_{max} > 0$. Moreover, for any given constant $\varphi \in \left(0, \frac{\pi}{2}\right)$, we define the set

$$S := \Re^n \times (-a, a)^n \times (-\varphi, \varphi)^n \times (0, v_{max})^n. \qquad (2.6)$$

The set $S$ in (2.6) represents all possible states of the system of all $n$ vehicles described by (2.3) and has the following interpretation. First, each vehicle should stay within the road, i.e., $(x_i, y_i) \in \Re \times (-a, a)$ for $i = 1, \ldots, n$. Moreover, with the given constant $\varphi \in \left(0, \frac{\pi}{2}\right)$, the vehicles should not be able to turn perpendicular to the road, as it should hold that $\theta_i \in (-\varphi, \varphi)$ for $i = 1, \ldots, n$. The constant $\varphi$ can be understood as a safety constraint, which allows the restriction of the movement of a vehicle; for instance, for vehicles moving at high speed, $\varphi$ should take values close to zero. Finally, the speeds of all vehicles should always be positive, i.e., no vehicle moves backwards at any time; and respect the road speed limits. One very important property, that is not captured by the set $S$, is that of collision avoidance between vehicles. This implies that the distance between the reference points of any pair of vehicles should always be greater than $L > 0$, where $L > 0$ is a safety distance that prevents collisions.

Due to the various constraints that were explained above, we must consider system (2.3) on the open set $\Omega \subset \Re^{4n}$ defined by:

$$\Omega := \{ w \in S : d_{i,j} > L, i, j = 1, \ldots, n, j \neq i \}. \qquad (2.7)$$

The set $\Omega$ in (2.7) describes the state-space of the $n$ vehicles operating on a lane-free road and acts as a basis for the problem formulation and for expressing the main objectives of the paper:

**Problem Statement**: For a group of $n$ vehicles modeled by (2.3) and operating on a lane-free road of width $2a > 0$, design decentralized feedback laws for $u_i$ and $F_i$ such that the following objectives hold:

(i) the vehicles do not collide with each other or with the boundary of the road, i.e., $d_{i,j}(t) > L$ for all $t \geq 0$, $i, j = 1, \ldots, n$, $j \neq i$, for a given constant $L > 0$, and $y_i(t) \in (-a, a)$ for all $t \geq 0$.
(ii) the speeds of all vehicles are always positive and remain below the given speed limit, i.e., $v_i(t) \in (0, v_{max})$ for $t \geq 0$, and converge to a given longitudinal speed set-point $v^* \in (0, v_{max})$, i.e., $\lim_{t \to +\infty} (v_i(t)) = v^*$, $i = 1, \ldots, n$.
(iii) the orientation of each vehicle is always bounded by the given value $\varphi \in \left(0, \frac{\pi}{2}\right)$, i.e., $\theta_i(t) \in (-\varphi, \varphi)$ for $t \geq 0$, and converges to zero, i.e., $\lim_{t \to +\infty} (\theta_i(t)) = 0$, $i = 1, \ldots, n$.
(iv) the accelerations, angular speeds, and lateral speeds of all vehicles tend to zero, i.e., $\lim_{t \to +\infty} (F_i(t)) = 0$, $\lim_{t \to +\infty} (u_i(t)) = 0$, and $\lim_{t \to +\infty} (\dot{y}_i(t)) = 0$, $i = 1, \ldots, n$.

It should be noted that, in mathematical terms, we require the closed-loop system to be well-posed on the state space $\Omega \subset \Re^{4n}$ defined by (2.7), i.e., for every initial condition $w(0) \in \Omega$, the closed-loop system (2.3), under the effect of all feedback laws for $u_i$ and $F_i$ for $i=1,...,n$, has a unique solution $w(t) \in \Omega$ defined for all $t \geq 0$. Moreover, we require that, for every initial condition $w(0) \in \Omega$, the solution $w(t) \in \Omega$ of the closed-loop system (2.3), under the effect of all feedback laws for $u_i$ and $F_i$ for $i=1,...,n$, satisfies $\lim_{t \to +\infty}(v_i(t)) = v^*$, $\lim_{t \to +\infty}(\theta_i(t)) = 0$, $\lim_{t \to +\infty}(F_i(t)) = 0$, $\lim_{t \to +\infty}(u_i(t)) = 0$ for all $i=1,...,n$. It should also be noticed that, since the lateral speed of each vehicle $\dot{y}_i$ is given by (2.3), it follows that the lateral speeds of all vehicles tend to zero, i.e., $\lim_{t \to +\infty}(\dot{y}_i(t)) = 0$ for $i=1,...,n$.

## 3. Main Results

3.1. Statements of main results

In this section, we design a novel decentralized control strategy in order to achieve the various objectives discussed in Section 2. First and foremost, we want to design the control inputs $u_i$ and $F_i$ in such a way that vehicles operating on a lane-free road do not collide with each other or with the boundary of the road. A typical approach for collision avoidance between vehicles is the use of repulsive potential functions (see for instance [5], [7], [8], [14] [23], [27], [29]). Repulsive potential functions are continuously differentiable functions, which repel vehicles based on their distance, with the force of repulsion being stronger as the distance between two vehicles becomes smaller, while there is little or no repulsion when the vehicles are distant. To that end, let $V:(L,+\infty) \to \mathfrak{R}_+$ be a $C^2$ function that satisfies:

$$\lim_{d \to L^+}(V(d)) = +\infty \tag{3.1}$$

$$V(d) = 0, \text{ for all } d \geq \lambda \tag{3.2}$$

where $\lambda > L$ is a constant. Let also $U:(-a,a) \to \mathfrak{R}_+$ be a $C^2$ function that satisfies:

$$\lim_{y \to (-a)^+}(U(y)) = +\infty, \quad \lim_{y \to a^-}(U(y)) = +\infty \tag{3.3}$$

$$U(0) = 0. \tag{3.4}$$

The potential function $U(y)$ in (3.3), (3.4) is designed so as to exert a repulsive force when the vehicles approach the boundary of the road.

To design feedback control laws that address objectives (i)-(iv) in the Problem Statement, we apply a control Lyapunov function methodology, where the feedback laws are selected appropriately to render the derivative of a Lyapunov function negative semi-definite. An appropriate function for this task is the following. Define, for all $w \in \Omega$,

$$H(w) := \frac{1}{2}\sum_{i=1}^{n}(v_i \cos(\theta_i) - v^*)^2 + \frac{1}{2}\sum_{i=1}^{n}v_i^2 \sin^2(\theta_i) + \sum_{i=1}^{n}U(y_i) + \frac{1}{2}\sum_{i=1}^{n}\sum_{j \neq i}V(d_{i,j}) + A\sum_{i=1}^{n}\left(\frac{1}{\cos(\theta_i) - \cos(\varphi)} - \frac{1}{1 - \cos(\varphi)}\right) \tag{3.5}$$

where $A>0$ is a parameter of the controller and the Lyapunov function, $v^* \in (0, v_{max})$ is a given longitudinal speed set-point, and $\varphi \in \left(0, \frac{\pi}{2}\right)$ is any constant that satisfies the inequality

$$\cos(\varphi) \geq \frac{v^*}{v_{max}}. \tag{3.6}$$

The function $H$ in (3.5), is inspired by the total energy of the system of $n$ vehicles and will allow us to exploit certain properties of the state space $\Omega$ in (2.7). The first two terms ($\frac{1}{2}\sum_{i=1}^{n}\left(v_i \cos(\theta_i) - v^*\right)^2 + \frac{1}{2}\sum_{i=1}^{n} v_i^2 \sin^2(\theta_i)$) represent the kinetic energy of the system of $n$ vehicles relative to an observer moving along the $x-$direction with speed equal to $v^*$; and they penalize the deviation of the longitudinal and lateral speeds from their desired values $v^*$ and zero, respectively. The sum of the third and fourth term ($\sum_{i=1}^{n} U(y_i) + \frac{1}{2}\sum_{i=1}^{n}\sum_{j \neq i} V(d_{i,j})$), which are based on the potential functions (3.3) and (3.4), is the potential energy of the system. Finally, the last term of (3.5) ($A\sum_{i=1}^{n}\left(\frac{1}{\cos(\theta_i)-\cos(\varphi)} - \frac{1}{1-\cos(\varphi)}\right)$) is a penalty term that blows up when $\theta_i \to \pm\varphi$. Inequality (3.6) is a technical assumption that restricts the movement of the vehicle when the desired speed is close to the road speed limit. Notice also that $H$ is not only a Lyapunov function, but possesses also certain characteristics of barrier functions, (see for instance [2], [12], [26]). Indeed, $H(w)$ grows unbounded on certain parts of the boundary of the state-space $\Omega$ in (2.7), i.e., when $y_i \to \pm a$ or $\theta_i \to \pm\varphi$ or $d_{i,j} \to L$ for some $i, j = 1,...,n$ with $i \neq j$ (recall (3.1) and (3.3)). The following result shows clearly that $H$ is a barrier function.

**Proposition 1:** *Let constants $A>0$, $v_{max}>0$, $v^* \in (0, v_{max})$, $\lambda > L > 0$, $\varphi \in \left(0, \frac{\pi}{2}\right)$ that satisfies (3.6), and define the function $H : \Omega \to \Re_+$ by means of (3.5), where $\Omega$ is given by (2.7). Then, there exist a non-decreasing function $\kappa : \Re_+ \to [0,a)$, a non-increasing function $\rho : \Re_+ \to (L, \lambda]$ and a non-decreasing function $\omega : \Re_+ \to [0, \varphi)$ such that the following implication holds:*

$$w \in \Omega \Rightarrow |\theta_i| \leq \omega(H(w)), |y_i| \leq \kappa(H(w)), d_{i,j} \geq \rho(H(w)), \text{ for } i, j=1,...,n, j \neq i. \tag{3.7}$$

The existence of a non-decreasing function $\kappa : \Re_+ \to [0,a)$, a non-increasing function $\rho : \Re_+ \to (L, \lambda]$ and a non-decreasing function $\omega : \Re_+ \to [0, \varphi)$, for which implication (3.7) holds, is a consequence of (3.1), (3.2), (3.3), (3.4) and definition (3.5). Implication (3.7) suggests that for any $w \in \Omega$, the orientations $\theta_i$ and the lateral positions $y_i$ of all vehicles $i=1,...,n$, as well as the distances $d_{i,j}$, $i, j=1,...,n, j \neq i$, are bounded by the energy of the system, see (3.5).

The feedback laws for each vehicle $i=1,...,n$ can be designed using (3.5), in terms of their own speed and orientation and the gradient of the potential functions $V_i$ and $U_i$ that satisfy (3.1), (3.2) and (3.3), (3.4), respectively:

$$u_i = -\left(v^* + \frac{A}{v_i\left(\cos(\theta_i) - \cos(\varphi)\right)^2}\right)^{-1}\left(\mu_1 v_i \sin(\theta_i) + U'(y_i) + p\sum_{j \neq i} V'(d_{i,j})\frac{(y_i - y_j)}{d_{i,j}} + \sin(\theta_i)F_i\right) \quad (3.8)$$

$$F_i = -\frac{k_i(w)}{\cos(\theta_i)}\left(v_i \cos(\theta_i) - v^*\right) - \frac{1}{\cos(\theta_i)}\sum_{j \neq i} V'(d_{i,j})\frac{(x_i - x_j)}{d_{i,j}} \quad (3.9)$$

$$k_i(w) = \mu_2 + \frac{1}{v^*}\sum_{j \neq i} V'(d_{i,j})\frac{(x_i - x_j)}{d_{i,j}} + \frac{v_{\max}\cos(\theta_i)}{v^*\left(v_{\max}\cos(\theta_i) - v^*\right)} f\left(-\sum_{j \neq i} V'(d_{i,j})\frac{(x_i - x_j)}{d_{i,j}}\right) \quad (3.10)$$

where $\mu_1, \mu_2 > 0$ are constants (controller gains) and $f \in C^1(\Re)$ is any function that satisfies

$$\max(x, 0) \leq f(x), \text{ for all } x \in \Re. \quad (3.11)$$

The term $k_i(w)$ in the acceleration $F_i(t)$, given by (3.9), is a state-dependent controller gain which guarantees that the speed of each vehicle will remain positive and less than the speed limit (see the proof of Theorem 1 below). The second term that appears in (3.9), is the summation of repelling forces ($V'(d)$) from vehicles that are in close proximity to vehicle $i$. If $V$ in (3.1), (3.2) is decreasing, then, the second term of (3.9) is positive if vehicle $j$ is behind vehicle $i$, i.e., $(x_i - x_j) > 0$. Indeed, in this case, we have that $-V'(d_{i,j})\frac{(x_i - x_j)}{d_{i,j}} > 0$, and this term represents the effect of nudging, since vehicles that are close and behind vehicle $i$ will exert a "pushing" force towards it that will increase its acceleration. It should be noticed that the control laws above are designed in such a way that the nudging force will not jeopardize traffic safety in terms of collisions, speeds exceeding desired bounds or vehicles departing from the road.

**Remark 1:**
   (i) Property (3.2) guarantees that the feedback laws (3.8), (3.9), (3.10) depend only on information from adjacent vehicles, namely from vehicles that are located at a distance less than $\lambda > 0$. Notice also that the control inputs (3.8), (3.9), (3.10) only require the distance from neighboring vehicles and not additional information, such as relative speeds ($v_i - v_j$) or relative orientations ($\theta_i - \theta_j$) $i, j = 1, ..., n, i \neq j$.
   (ii) Any function $f \in C^1(\Re)$ that satisfies (3.11) can be used in (3.10). For example, the function $f(x) = \frac{\varepsilon}{2} + \frac{1}{2\varepsilon}x^2$ for every $\varepsilon > 0$ satisfies (3.11), since $\max(x, 0) \leq |x| \leq \frac{\varepsilon}{2} + \frac{1}{2\varepsilon}x^2$ for all $x \in \Re$. Another function that satisfies (3.11) is the function

$$f(x) = \frac{1}{2\varepsilon}\begin{cases} 0 & \text{if } x \leq -\varepsilon \\ (x + \varepsilon)^2 & \text{if } -\varepsilon < x < 0 \\ \varepsilon^2 + 2\varepsilon x & \text{if } x \geq 0 \end{cases} \quad (3.12)$$

for every $\varepsilon > 0$. This generic design for the function $f$ will allow to regulate the longitudinal acceleration as desired. For instance, in the first example above, $f$ exhibits quadratic growth while in (3.12) only linear growth for $x \geq 0$.

Let $p \geq 1$, and consider two concentric ellipses with semi-major axes $L$ and $\lambda$, with $L < \lambda$, and semi-minor axes $\frac{L}{\sqrt{p}}$ and $\frac{\lambda}{\sqrt{p}}$, respectively. Let $m \geq 2$ be the maximum number of points that can be placed within the area bounded by the two concentric ellipses, so that each point has distance (in the metric given by (2.4)) at least $L$ from every other point. The following proposition presents certain properties of the control laws (3.8), (3.9), (3.10) and their relation with the speed limit $v_{max}$, the longitudinal speed set-point $v^*$, and the maximum number of neighboring vehicles $m$.

**Proposition 2:** *Let constants $\lambda > L > 0$, $a > 0$, $p \geq 1$, and let $V : (L, +\infty) \to \Re_+$, $U : (-a, a) \to \Re_+$ be $C^2$ functions that satisfy (3.1), (3.2) and (3.3), (3.4), respectively, and define*

$$b_1(s) := \max\{|V'(d)| : s \leq d \leq \lambda\} \text{ for } s \in (L, \lambda] \tag{3.13}$$

$$b_2(s) := \max\{|U'(y)| : |y| \leq s\} \text{ for } s \in [0, a). \tag{3.14}$$

*Define the set $\Omega$ by means of (2.7). Then, for any $w \in \Omega$, there exist a non-decreasing function $\kappa : \Re_+ \to [0, a)$ and a non-increasing function $\rho : \Re_+ \to (L, \lambda]$ such that the functions $u_i$, $F_i$ and $k_i$ in (3.8), (3.9), and (3.10), respectively, satisfy the following inequalities*

$$k_i(w)v^* \geq \sum_{j \neq i} V'(d_{i,j}) \frac{(x_i - x_j)}{d_{i,j}} \geq -k_i(w)\left(v_{max} \cos(\theta_i) - v^*\right), \, i = 1, \ldots, n; \tag{3.15}$$

$$\mu_2 \leq k_i(w) \leq R(H(w)), \, i = 1, \ldots, n; \tag{3.16}$$

$$k_i(w)v_{max} \geq k_i(w)(v_{max} - v_i) \geq F_i \geq -k_i(w)v_i \geq -k_i(w)v_{max}; \tag{3.17}$$

$$|u_i| \leq \frac{1}{v^*}\left((\mu_1 + k_i(w))v_{max} + b_2(\kappa(H(w))) + m\sqrt{p}\, b_1(\rho(H(w)))\right) \tag{3.18}$$

*where $R : \Re_+ \to \Re_+$ is the increasing function defined by*

$$R(s) := \mu_2 + \frac{m}{v^*}b_1(\rho(s)) + v_{max} \frac{(A + s\cos(\varphi)(1 - \cos(\varphi)))\max\{f(z) : |z| \leq mb_1(\rho(s))\}}{Av^*(v_{max} - v^*) + v^*(v_{max}\cos(\varphi) - v^*)(1 - \cos(\varphi))s}. \tag{3.19}$$

The proof of inequalities (3.15), (3.16), (3.17), and (3.18) follows by exploiting (3.5), (3.8), (3.9), (3.10), and the result of Proposition 1. Inequality (3.16) suggests that the magnitude of $k_i(w)$ depends on the "energy" of the system defined by the Lyapunov-like function $H$ and the maximum number of neighboring vehicles $m \geq 2$. Moreover, $k_i(w)$ plays an important role, since it provides certain bounds on the acceleration $F_i$ in (3.17) and the maximum "nudging" effect that each vehicle $i = 1, \ldots, n$ experiences from neighboring vehicles as described in (3.15).

We are now in a position to state the main result.

**Theorem 1:** *Suppose that there exist constants $a > 0$, $\lambda > L > 0$, $p \geq 1$ and $C^2$ functions $V : (L, +\infty) \to \Re_+$, $U : (-a, a) \to \Re_+$ that satisfy (3.1), (3.2), and (3.3), (3.4), respectively. In addition, for given constants $v_{\max} > 0$, $v^* \in (0, v_{\max})$, and $\varphi \in \left(0, \dfrac{\pi}{2}\right)$ that satisfies (3.6), define the function $H : \Omega \to \Re_+$ by means of (3.5) where $\Omega$ is given by (2.7). Then, for every $w_0 \in \Omega$ there exists a unique solution $w(t) \in \Omega$ of the initial-value problem (2.3), (3.8), (3.9), (3.10) with initial condition $w(0) = w_0$. The solution $w(t) \in \Omega$ is defined for all $t \geq 0$ and satisfies for $i = 1, \ldots, n$*

$$\lim_{t \to +\infty}(v_i(t)) = v^*, \quad \lim_{t \to +\infty}(\theta_i(t)) = 0 \tag{3.20}$$

$$\lim_{t \to +\infty}(u_i(t)) = 0, \quad \lim_{t \to +\infty}(F_i(t)) = 0. \tag{3.21}$$

*Moreover, there exist a non-decreasing function $\kappa : \Re_+ \to [0, a)$ and a non-increasing function $\rho : \Re_+ \to (L, \lambda]$ such that*

$$|F_i(t)| \leq R(H(w_0))v_{\max}, \text{ for all } t \geq 0 \tag{3.22}$$

$$|u_i(t)| \leq \frac{1}{v^*}\left((\mu_1 + R(H(w_0)))v_{\max} + b_2(\kappa(H(w_0))) + m\sqrt{p}b_1(\rho(H(w_0)))\right), \text{ for all } t \geq 0 \tag{3.23}$$

*where $b_1, b_2, R$ are defined by (3.13), (3.14) and (3.19), respectively.*

**Remark 2:**
(i) It is important to notice that due to technical constraints, an inequality of the form $|F_i(t)| \leq K$ must be satisfied for all $t \geq 0$, where $K > 0$ is a constant that depends on the technical characteristics of the vehicles and the road. Inequality (3.22) allows us to determine the set of initial conditions $w_0 \in \Omega$ for which the inequality $|F_i(t)| \leq K$ holds: it includes the set of all $w_0 \in \Omega$ with $R(H(w_0))v_{\max} \leq K$.

(ii) Although we cannot predict the "ultimate" arrangement of the vehicles on the road (and we cannot even show that a final configuration of the vehicles on the road is attained; see remark below), the limits (3.20), (3.21) and definitions (3.9) allow us to predict that

$$\lim_{t \to +\infty}\left(\sum_{j \neq i} V'(d_{i,j}(t))\frac{(x_i(t) - x_j(t))}{d_{i,j}(t)}\right) = \lim_{t \to +\infty}\left(U'(y_i(t)) + p\sum_{j \neq i}V'(d_{i,j}(t))\frac{(y_i(t) - y_j(t))}{d_{i,j}(t)}\right) = 0 \quad \text{for}$$

$i = 1, \ldots, n$. Consequently, the "ultimate" arrangement of the vehicles in the road (if such a thing exists) must satisfy the equations $\sum_{j \neq i} V'(d_{i,j})\dfrac{(x_i - x_j)}{d_{i,j}} = U'(y_i) + p\sum_{j \neq i}V'(d_{i,j})\dfrac{(y_i - y_j)}{d_{i,j}} = 0$ for $i = 1, \ldots, n$ as well as the constraints $|y_i| < a$, $d_{i,j} > L$ for $i, j = 1, \ldots, n$, $j \neq i$. Despite the fact that the constrained system of $2n$ equations has infinite solutions, not every arrangement of vehicles satisfies the aforementioned constrained system.

(iii) The proof of Theorem 1 relies on Barbălat's lemma [13] and does not use LaSalle's invariance principle. The reason that LaSalle's invariance principle cannot be used for the proof of Theorem 1 is the fact that the state components $x_i(t)$, $i = 1, \ldots, n$, do not take values in a bounded set. Moreover, we cannot show that the relative positions of the vehicles, i.e., the quantities $x_i(t) - x_j(t)$ for $i, j = 1, \ldots, n$, $j \neq i$, take values in a bounded set. Thus, we cannot show that the limits $\lim_{t\to+\infty}(y_i(t))$, $\lim_{t\to+\infty}(x_i(t)-x_j(t))$, $\lim_{t\to+\infty}(d_{i,j}(t))$ for $i,j=1,...,n$, $j\neq i$, exist. Consequently, we cannot ensure that a final configuration of the vehicles on the road will be attained. However, the proof of Theorem 1 shows that $\lim_{t\to+\infty}(\dot{y}_i(t))=0$, $\lim_{t\to+\infty}(\dot{x}_i(t)-\dot{x}_j(t))=0$, $\lim_{t\to+\infty}(\dot{d}_{i,j}(t))=0$ for $i,j=1,...,n$, $j\neq i$ (a consequence of (2.3), (3.20) and inequality (5.29) in the proof of Theorem 1). Therefore, it is expected that the vehicles on the road will approach a final configuration.

### 3.2. Optimal Selection of Safety Distance and Eccentricity

One of the main advantages of a lane-free traffic is the improved exploitation of the width of the road. Since the distance between vehicles is measured according to (2.4), the number of vehicles $N$ that can be placed side-by-side is directly related to the safety distance $L>0$ and the weight $p\geq 1$ in (2.4) according to the formula

$$N = \frac{2a\sqrt{p}}{L} \tag{3.24}$$

which is obtained by dividing the width of the road $2a>0$ by $\frac{L}{\sqrt{p}}$ which is the semi-minor axis of the ellipse $x_i^2 + py_i^2 = L^2$ placed around each vehicle $i=1,...,n$. We discuss next the optimal selection of the safety distance $L>0$ and of the weight factor $p\geq 1$ for a given bound $\varphi \in \left(0, \frac{\pi}{2}\right)$ on the orientation of the vehicles that satisfies (3.6).

<u>Selection of safety distance:</u> Let $p\geq 1$ and $\varphi \in \left(0, \frac{\pi}{2}\right)$. Then, the maximum distance, at which two identical vehicles of length $\sigma > 0$ may collide, is

$$L = \sigma \max\left(2\sqrt{p}\sin(\varphi), \sqrt{1+(p-1)\sin^2(\varphi)}\right). \tag{3.25}$$

Finding the maximum distance at which two vehicles collide, given their length and maximum orientation angle, is a geometric problem and its solution is based on the intersection of the line segments $I_1 = \{(\eta_1\sigma\cos(\theta_1), \eta_1\sigma\sin(\theta_1)): \eta_1 \in [0,1]\}$ and $I_2 = \{(x+\eta_2\sigma\cos(\theta_2), y+\eta_2\sigma\sin(\theta_2)): \eta_2 \in [0,1]\}$ with $|\theta_i|\leq \varphi < \frac{\pi}{2}$, $i=1,2$. Intuitively, each line segment represents a vehicle of length $\sigma > 0$ that is placed at the positions $(0,0)$ and $(x,y)$, respectively. Hence, given their orientations $|\theta_i|\leq \varphi < \frac{\pi}{2}$, $i=1,2$, a non-empty intersection of the segments $I_1$ and $I_2$, indicates a collision between the two vehicles. Notice that $I_1 \cap I_2 \neq \emptyset$, if there exist $\eta_1, \eta_2 \in [0,1]$ with $\eta_1\sigma\cos(\theta_1) = x+\eta_2\sigma\cos(\theta_2)$, $\eta_1\sigma\sin(\theta_1) = y+\eta_2\sigma\sin(\theta_2)$. Then, the distance, as defined in (2.4), between the points $(0,0)$ and $(x,y)$ is $d = \sqrt{x^2+py^2} = \sigma\sqrt{(\eta_1\cos(\theta_1)-\eta_2\cos(\theta_2))^2 + p(\eta_1\sin(\theta_1)-\eta_2\sin(\theta_2))^2}$. Condition (3.25) can be obtained by using the following lemma to maximize the distance $d$.

**Lemma 1:** Let $p \geq 1$ and $\varphi \in \left(0, \frac{\pi}{2}\right]$. Let also, $\theta_1, \theta_2 \in [-\varphi, \varphi]$ and define

$$f(\eta_1, \eta_2) = (\eta_1 \cos(\theta_1) - \eta_2 \cos(\theta_2))^2 + p(\eta_1 \sin(\theta_1) - \eta_2 \sin(\theta_2))^2. \qquad (3.26)$$

Then, for any $\eta_1, \eta_2 \in [0,1]$, it holds that

$$f(\eta_1, \eta_2) \leq \max\left(4p\sin^2(\varphi), 1 + (p-1)\sin^2(\varphi)\right). \qquad (3.27)$$

Formula (3.25) provides explicitly the safety distance $L$ between vehicles for given constant $p \geq 1$ and the maximum orientation angle $\varphi \in \left(0, \frac{\pi}{2}\right]$. Since the safety distance defines an ellipse $x_i^2 + py_i^2 = L^2$ around every vehicle $i = 1, \ldots, n$, the eccentricity $e = \sqrt{1 - \frac{1}{p}}$ for each ellipse can be determined by appropriately selecting the constant $p \geq 1$. In order to maximize the number of vehicles that can be placed side-by-side $N$ in (3.24), it suffices to minimize the semi-minor axis $\frac{L}{\sqrt{p}}$ of the ellipse. It follows using (3.25) that the optimal selection of $p$ is

$$\min_{p \geq 1}\left(\frac{L}{\sqrt{p}}\right) = 2\sigma \sin(\varphi), \text{ attained for } p = \frac{1}{3\tan^2(\varphi)} \text{ when } \varphi \leq \frac{\pi}{6} \text{ and } p = 1 \text{ when } \varphi > \frac{\pi}{6}.$$

Notice that for given $\varphi \leq \frac{\pi}{6}$ (or $\varphi > \frac{\pi}{6}$), larger values of $p > \frac{1}{3\tan^2(\varphi)}$ (or $p > 1$, respectively) will not affect the minimum value of $\frac{L}{\sqrt{p}}$, see Figure 2.

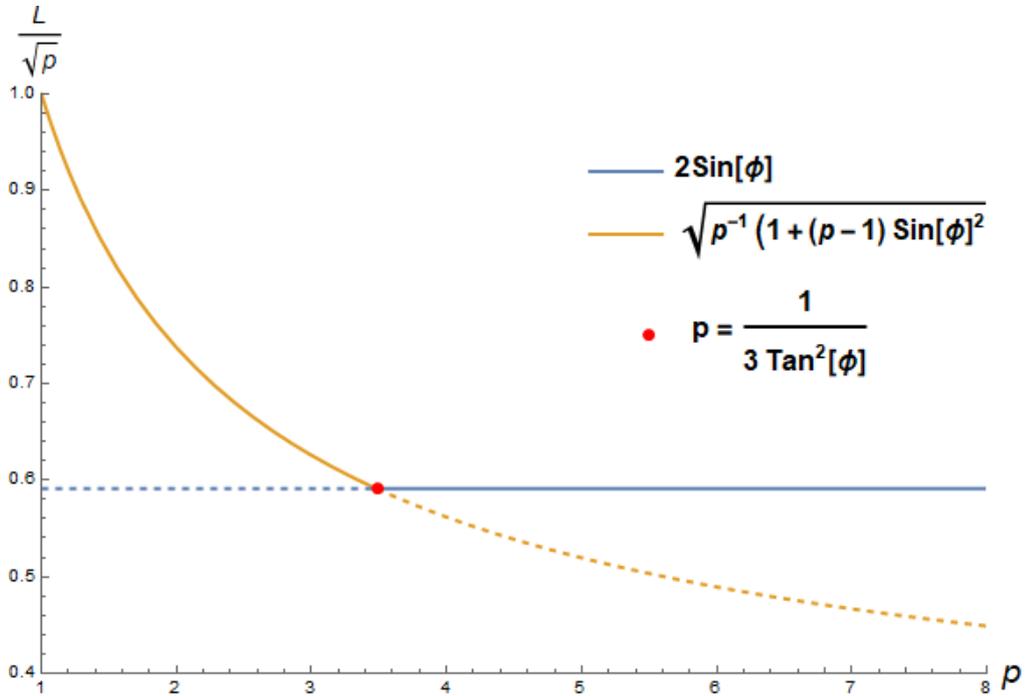

**Figure 2:** Optimal selection of eccentricity.

# 4. Illustrative Examples

In the simulation results below, we demonstrate the application and effectiveness of the proposed decentralized cruise controllers for autonomous vehicles driving on lane-free roads. Specifically, we consider a group of $n=10$ vehicles on a lane-free road of width $2a>0$, modeled as in (2.2) with the feedback laws (3.8), (3.9), (3.10), and $f(x)$ given by means of (3.12). The vehicle-repulsive potential function $V$ and the boundary-repulsive potential function $U$ are specified as

$$V(d) = \begin{cases} q\frac{(\lambda-d)^3}{d-L} &, L < d \leq \lambda \\ 0 &, d > \lambda \end{cases},$$

$$U(y) = \begin{cases} \left(\frac{1}{a^2-y^2} - \frac{c}{a^2}\right)^4 &, -a < y < -\frac{a\sqrt{c-1}}{\sqrt{c}} \text{ and } \frac{a\sqrt{c-1}}{\sqrt{c}} < y < a \\ 0 &, -\frac{a\sqrt{c-1}}{\sqrt{c}} \leq y \leq \frac{a\sqrt{c-1}}{\sqrt{c}} \end{cases}$$

where $c \geq 1, q > 0$ are design parameters. Notice that $V$ and $U$ above, satisfy (3.1), (3.2) and (3.3), (3.4), respectively. By appropriately selecting the constant $q$, we can adjust the repulsion force of the potential $V$ and consequently the magnitude of the acceleration $F_i$, see (3.9). In particular, for small values of $q$, the values of $V$ (and consequently the acceleration $F_i$) will be smaller away from $L$, but will increase more sharply as $d \to L$. The constant $c \geq 1$ affects the final configuration of the vehicles relative to the boundary of the road. More specifically, for $c=1$ we have that $U(y)=0$ if $y=0$, which will force the vehicles to form a single platoon in the middle of the road. For $c>1$, we have that $U(y)=0$ in a neighborhood around $y=0$, and the vehicles' converged lateral positions in this case will be within the strip $-\frac{a\sqrt{c-1}}{\sqrt{c}} \leq y \leq \frac{a\sqrt{c-1}}{\sqrt{c}}$.

To verify numerically and illustrate the results of Theorem 1, we assume that all vehicles have length $\sigma = 5m$ and operate on a road with speed limit $v_{\max} = 35m/s$ and width $2a = 14.4m$, which corresponds to a road with 4 conventional lanes of width $3.6m$. We set the longitudinal set-point $v^* = 30m/s$ and select $\varphi = 0.25$ in order to satisfy condition (3.6). Using (3.18) and (3.25), we obtain the optimal eccentricity and safety distance $p = 5.11$ and $L = 5.59m$, respectively. This choice allows us to effectively use the full width of the road and increases the lateral occupancy by $45\%$, since we get from (3.24) that $N = 5.8$. We set $\varepsilon = 0.2$, $\mu_1 = 0.5$, $\mu_2 = 0.1$, $q = 3*10^{-3}$, $\lambda = 25m$, $A=1$, and $c = 1.5$. Assigning small values of the constants $\varepsilon$ and $\mu_2$ in (3.12) and (3.10), respectively, provide smaller values for the acceleration $F_i$. Analogously, small values of the constants $\mu_1$ and $A$ in (3.8) lead to a smooth rotation rate, without abrupt turns. Each vehicle is assigned an index $i=1,\ldots,10$ and a unique color (as shown at the top of Figure 3), which will be used in every figure that appears hereafter.

Figure 3 shows snapshots of the simulation, with the vehicles' positions and orientations being visible at different time instants. The small ellipse around each vehicle is the safety region defined by the ellipse $x^2 + py^2 = L^2$, and the larger ellipse corresponds to the area defined by the ellipse $x^2 + py^2 = \lambda^2$. It is seen that all vehicles remain within the road and do not collide with each other.

The final snapshot at $t = 340s$ shows that the vehicles move towards a configuration with lateral positions $-4.15 \leq y \leq 4.15$.

Figure 4 displays the longitudinal speed $\dot{x}_i$ and the acceleration $F_i$ of each vehicle. The speeds of all vehicles remain within the bounds $(0, v_{max})$ and converge to the longitudinal set-point $v^*$. In Figure 5, we focus only on the vehicle 7 (cyan) and vehicle 9 (purple). It can be seen from Figure 3 and Figure 5, that vehicle 7 approaches vehicle 9 at much higher speed. Notice that the speed and acceleration of vehicle 9 is increased, as the distance from the vehicle 7 decreases. This is exactly the effect of nudging, i.e., vehicle 7 exerts a pushing force on vehicle 9, which, as a result, increases its acceleration and speed. On the other hand, it can be seen in Figure 5, that the vehicle 7 decelerates (is repulsed) to avoid collision with vehicle 9.

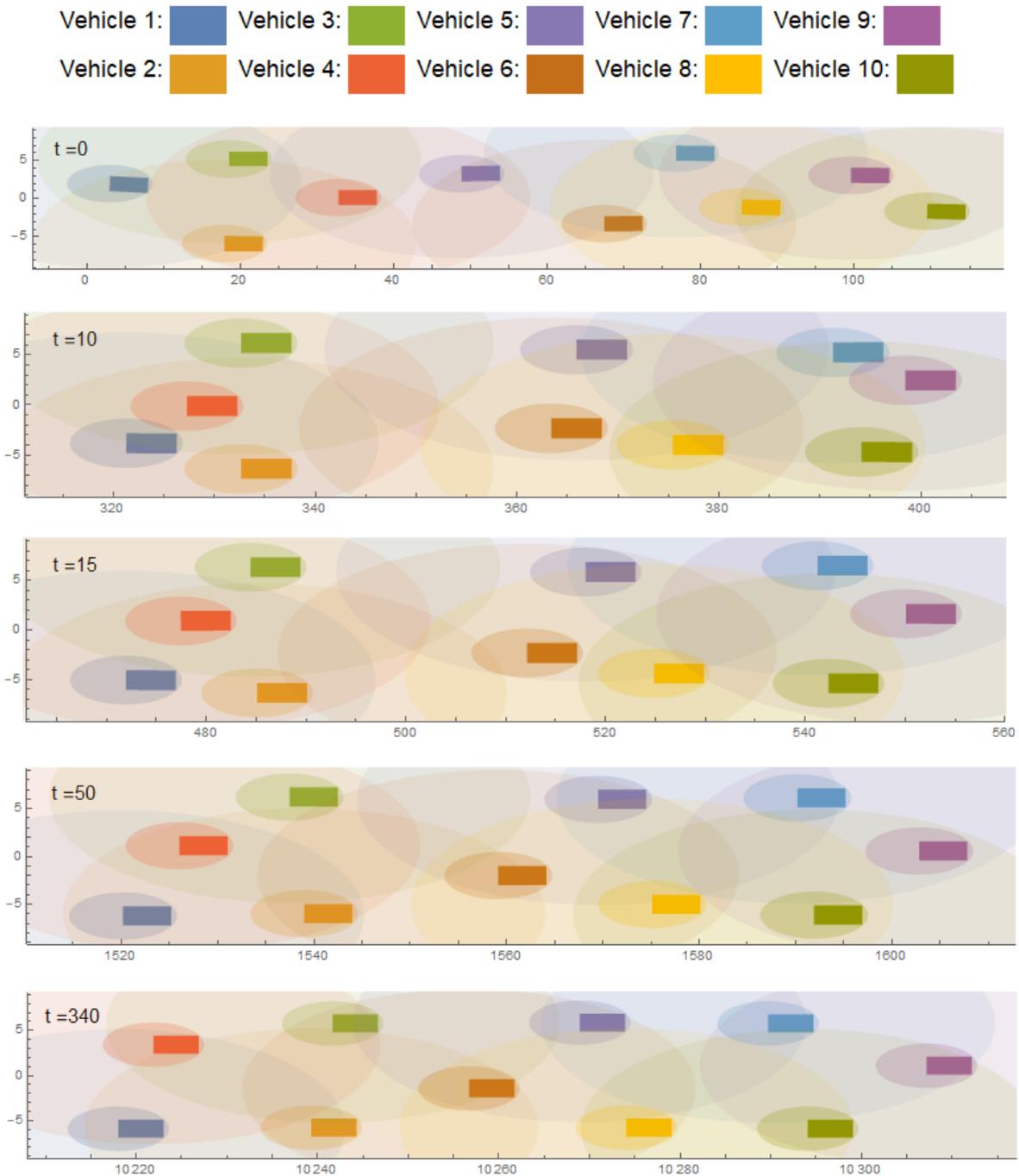

**Figure 3:** Snapshots at various time instants of 10 vehicles operating on a lane-free road.

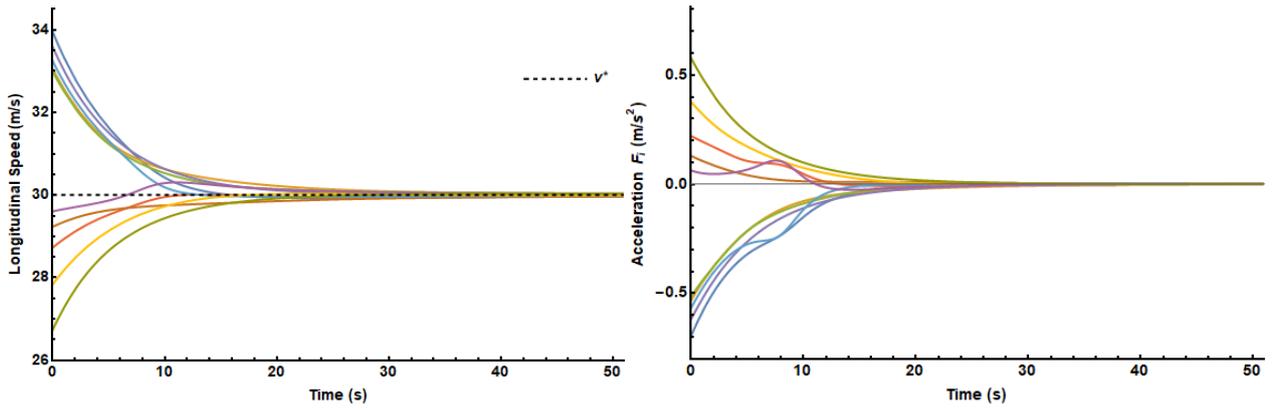

**Figure 4:** The longitudinal speed of each vehicle where all speeds converge to the speed set-point $v^*$ on the left; and the acceleration $F_i$ of each vehicle on the right.

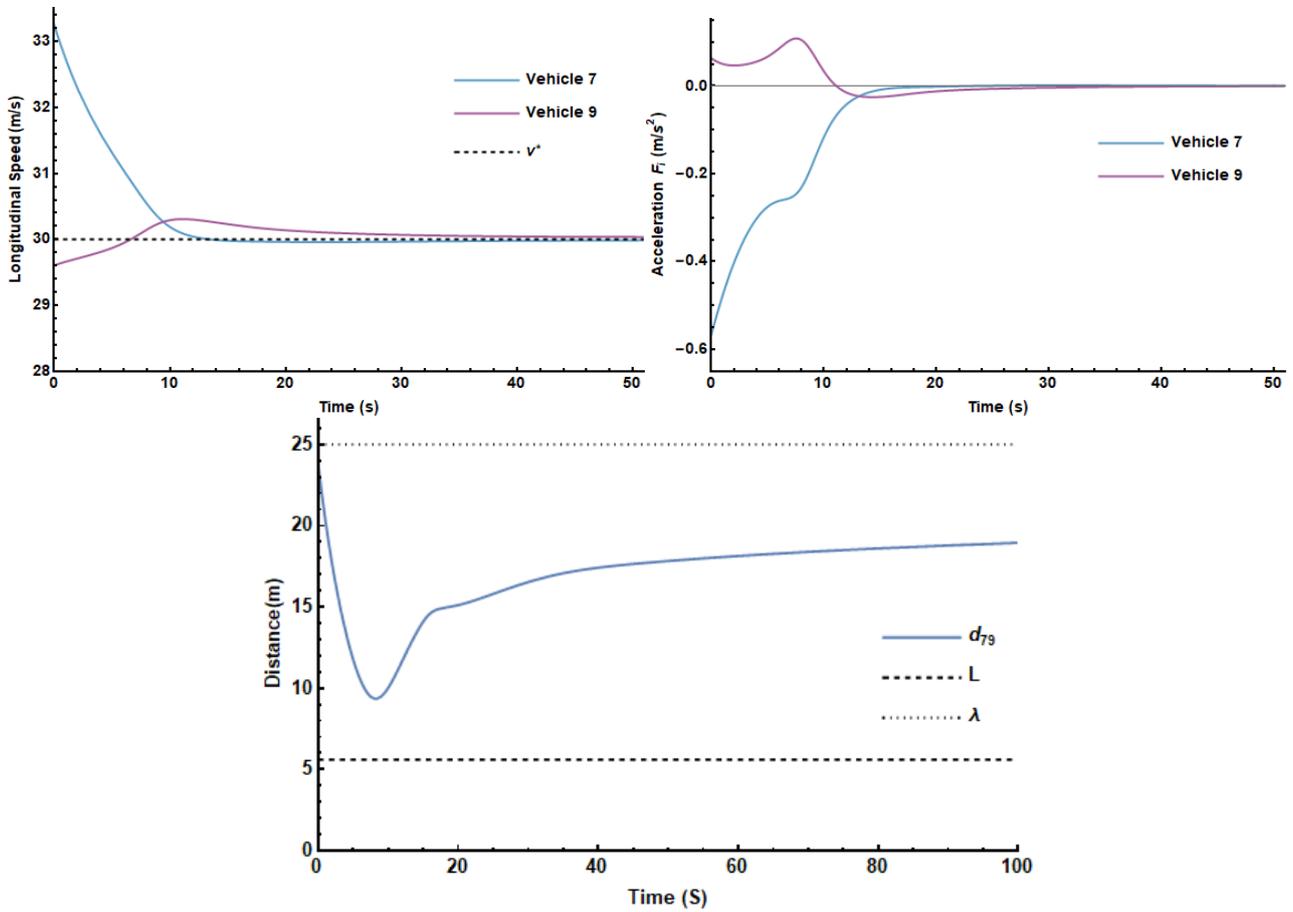

**Figure 5:** The effect of nudging. Vehicle 9 accelerates as the (elliptical) distance from vehicle 7 decreases.

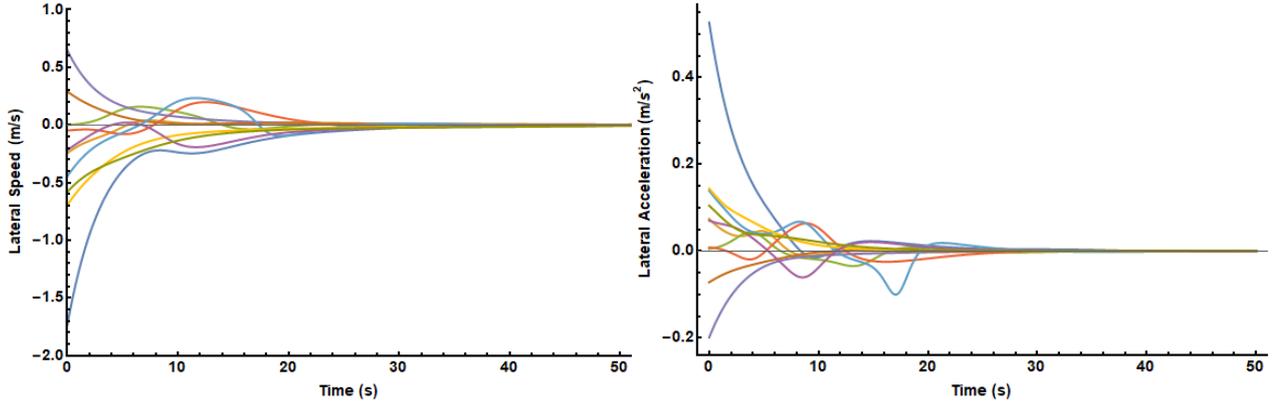

**Figure 6:** The lateral speed (left) and lateral acceleration (right) of each vehicle.

The lateral speeds $\dot{y}_i$ and lateral accelerations $\ddot{y}_i$ of the vehicles are shown in Figure 6; both converge to zero, indicating that eventually the vehicles move parallel to the road. The observed temporarily changing values of some lateral accelerations and speeds near the times $t = 8s$ and $t = 18s$ occur because of interactions when the vehicles approach the boundary of the road or approach other vehicles. For example, the temporary increase of the lateral acceleration of vehicle 7 and simultaneous temporary decrease of lateral acceleration of vehicle 9, around time $t = 8s$, are an aftereffect of vehicle interaction, whereby the vehicle 7 is nudged and accelerates towards the right side of the road, while the nudging vehicle 9 is repulsed. An analogous interaction can also be observed for vehicle 1, vehicle 2, and vehicle 4, at the time instant t=8s, see Figure 3 and Figure 6.

Figure 7 shows the rotation rates $u_i$ and the orientations $\theta_i$, all converging to zero as suggested by Theorem 1. Again, the vehicles are seen to change their rotation rate and orientation, which lead to the aforementioned changes in lateral speed and lateral acceleration, when they approach the boundary of the road or other vehicles.

Finally, Figure 8 depicts the minimum inter-vehicle distance $d_{\min}(t) := \min\{d_{i,j}(t), i, j = 1,...,n, i \neq j\}$ (blue line), showing that the vehicles do not collide with each other, since $d_{i,j}(t) > L$, $i, j = 1,...,n, i \neq j$ at all times. Moreover, Figure 8 also shows the minimum inter-vehicle distance using the same initial conditions with $\lambda = 40m$ (yellow line).

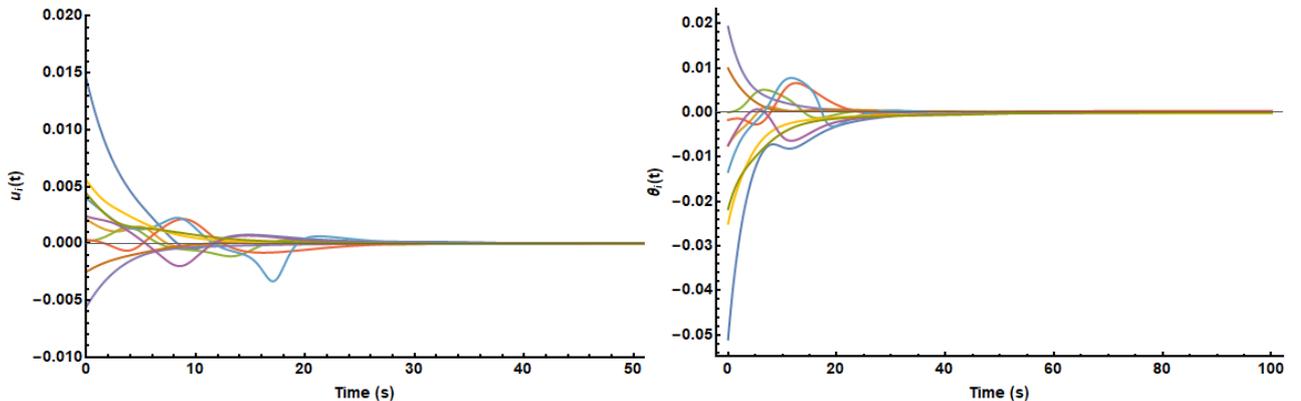

**Figure 7:** The rotation rate $u_i$ and orientation $\theta_i$ converge to zero as indicated by Theorem 1.

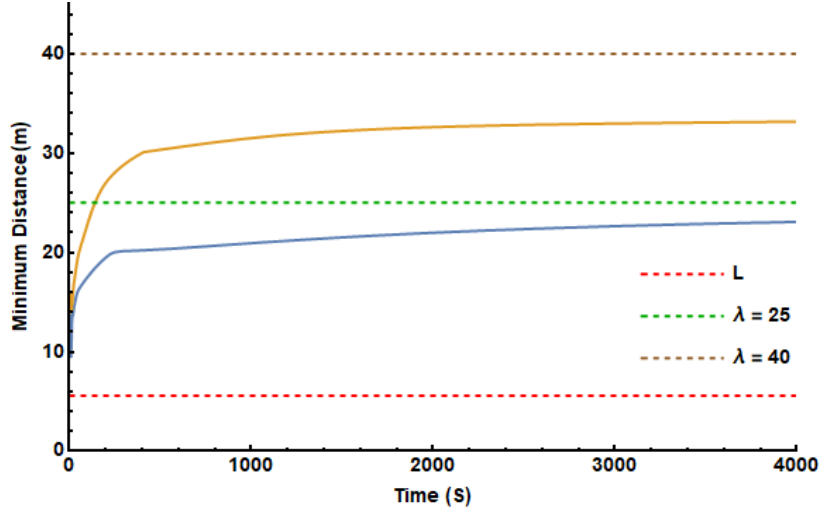

**Figure 8:** The minimum inter-vehicle distance with blue for $\lambda = 25m$ and with yellow for $\lambda = 40m$, which verify that there are no collisions among vehicles.

Figure 8 is representative for the controller behavior, as it was found, in multiple simulations with different values of $\lambda$, that the effect of increasing the value of $\lambda$ is an increase of the minimum inter-vehicle distance $d_{min}$ over time, as well as in the converged state.

## 5. Proofs of Main Results

For the proof of Proposition 1, the following technical lemma is required.

**Lemma 2:** *Let $c > 0$ be a given constant and let $\beta : [0,c) \to \Re_+$ be a non-decreasing function with $\beta(0) = 0$ and $\lim_{s \to c^-}(\beta(s)) = +\infty$. Then there exists a constant $s^* \in [0,c)$ and a continuous, non-decreasing function $\bar{\beta} : [0,c) \to \Re_+$, for which the following properties hold:*
  *a) $\lim_{s \to c^-}(\bar{\beta}(s)) = +\infty$, $\bar{\beta}(s^*) = 0$,*
  *b) $\bar{\beta}(s) \leq \beta(s)$ for all $s \geq 0$,*
  *c) $\bar{\beta}$ is increasing on $[s^*, c)$.*

**Proof of Lemma 2:** We first extend the domain of $\beta : [0,c) \to \Re_+$ by setting $\beta(s) = 0$ for $s < 0$. We define for $s \in [0,c)$:

$$\bar{\beta}(s) := \exp(s-c) \begin{cases} c \int_{c\frac{2s-c}{s}}^{s} \frac{\beta(u)\,du}{(c-u)^2} & \text{if } s \in (0,c) \\ 0 & \text{if } s = 0 \end{cases} \quad (5.1)$$

Notice that monotonicity of $\beta$ and definition (5.1) guarantee that $\bar{\beta} : [0,c) \to \Re_+$ is non-decreasing and satisfies the following estimates:

$$\exp(s-c)\beta\left(c\frac{2s-c}{s}\right) \leq \bar{\beta}(s) \leq \exp(s-c)\beta(s) \text{ for } s > 0 \quad (5.2)$$

$$\bar{\beta}(s) \leq \exp(s-c) \frac{s\beta(s)}{a-s}, \text{ for } s \in \left[0, \frac{c}{2}\right]. \tag{5.3}$$

Definition (5.1) and estimates (5.2), (5.3) guarantee that $\bar{\beta}:[0,c) \to \Re_+$ is continuous with $\lim_{s \to c^-}(\bar{\beta}(s)) = +\infty$. The inequality $\bar{\beta}(s) \leq \beta(s)$ for all $s \geq 0$ is a direct consequence of (5.2) and definition (5.1). Finally, define $s^* := \sup\{r \in [0,c): \bar{\beta}(r) = 0\}$. Since $\lim_{s \to c^-}(\bar{\beta}(s)) = +\infty$, it follows that $0 \leq s^* < c$. Definition (5.1) guarantees that $\bar{\beta}:[s^*,c) \to \Re_+$ is increasing. Moreover, continuity of $\bar{\beta}$ implies that $\bar{\beta}(s^*) = 0$. The proof is complete. ◁

**Proof of Proposition 1:** Using (3.5), we get for all $w \in \Omega$ and $i, j = 1,...,n$ with $i \neq j$:

$$A\left(\frac{1}{\cos(\theta_i) - \cos(\varphi)} - \frac{1}{1 - \cos(\varphi)}\right) \leq H(w) \tag{5.4}$$

$$U(y_i) \leq H(w) \tag{5.5}$$

$$V(d_{i,j}) \leq H(w). \tag{5.6}$$

Inequality (5.4) and equation (3.6) imply that $|\theta_i| \leq \arccos\left(\cos(\varphi) + \frac{A(1-\cos(\varphi))}{A + (1-\cos(\varphi))H(w)}\right) < \varphi$ and consequently, (3.7) holds with $\omega(s) := \arccos\left(\cos(\varphi) + \frac{A(1-\cos(\varphi))}{A + (1-\cos(\varphi))s}\right)$.

Define
$$\tilde{b}(s) := \min\{U(y): y \in (-a,-s] \cup [s,a)\}, \text{ for } s \in [0,a). \tag{5.7}$$

Notice that (3.3), (3.4), and definition (5.7) imply that $\tilde{b}:[0,a) \to \Re_+$ is non-decreasing with $\tilde{b}(0) = 0$ and $\lim_{s \to a^-}(\tilde{b}(s)) = +\infty$. We apply Lemma 2 with $\beta = \tilde{b}$, $c = a$ and we conclude that there exists a constant $s^* \in [0,a)$ and a continuous, non-decreasing function $\bar{\beta}:[0,a) \to \Re_+$, with $\lim_{s \to a^-}(\bar{\beta}(s)) = +\infty$, $\bar{\beta}(s^*) = 0$, $\bar{\beta}(s) \leq \tilde{b}(s)$ for all $s \geq 0$ and $\bar{\beta}$ being increasing on $[s^*,a)$. Consider the inverse function $\kappa: \Re_+ \to [s^*,a)$ of $\bar{\beta}$ restricted on $[s^*,a)$. Notice that the inverse function exists because $\bar{\beta}$ is increasing on $[s^*,a)$ with $\lim_{s \to a^-}(\bar{\beta}(s)) = +\infty$, $\bar{\beta}(s^*) = 0$, it is continuous function because $\bar{\beta}$ is continuous and is increasing because $\bar{\beta}$ is increasing. It follows from (5.5) and definition (5.7) that $\bar{\beta}(|y_i|) \leq \tilde{b}(|y_i|) \leq H(w)$. The previous inequality shows that $|y_i| \leq \kappa(H(w)) < a$.

Next, define
$$\theta(s) := \min\{V(d): L < d \leq s\}, \text{ for } s > L. \tag{5.8}$$

Notice that (3.1), (3.2) and definition (5.8) imply that $\theta:(L,+\infty)\to\mathfrak{R}_+$ is non-increasing with $\theta(s)=0$ for $s\geq\lambda$ and $\lim_{s\to L^+}(\theta(s))=+\infty$. It follows that the function $\tilde{\theta}:[0,1/L)\to\mathfrak{R}_+$, defined by $\tilde{\theta}(s)=\theta(1/s)$ for $s>0$ and $\tilde{\theta}(0)=0$, is a non-decreasing function with $\lim_{s\to(1/L)^-}(\tilde{\theta}(s))=+\infty$. Moreover, if holds that $\tilde{\theta}(s)=0$ for $s\in[0,1/\lambda]$. We apply again Lemma 2 with $\beta=\tilde{\theta}$, $c=1/L$ and we conclude that there exists a constant $s^{**}\in[0,1/L)$ and a continuous, non-decreasing function $\bar{\theta}:[0,1/L)\to\mathfrak{R}_+$, with $\lim_{s\to(1/L)^-}(\bar{\theta}(s))=+\infty$, $\bar{\theta}(s^{**})=0$, $\bar{\theta}(s)\leq\tilde{\theta}(s)$ for all $s\geq 0$ and $\bar{\theta}$ being increasing on $[s^{**},1/L)$. Since $\bar{\theta}(s)>0$ for all $s>s^{**}$, the constant $s^{**}$ must be greater than or equal to $1/\lambda$. Consider the inverse function $\gamma:\mathfrak{R}_+\to[s^{**},1/L)$ of $\bar{\theta}$ restricted on $[s^{**},1/L)$. Notice that the inverse function exists because $\bar{\theta}$ is increasing on $[s^{**},1/L)$ with $\lim_{s\to(1/L)^-}(\bar{\theta}(s))=+\infty$, $\bar{\theta}(s^{**})=0$, it is continuous function because $\bar{\theta}$ is continuous and is increasing because $\bar{\theta}$ is increasing.

It follows from (5.6) and definition (5.8) that $\theta(d_{i,j})\leq V(d_{i,j})\leq H(w)$. Since $\tilde{\theta}(s)=\theta(1/s)$ for $s>0$, the previous inequality implies the inequality $\tilde{\theta}(1/d_{i,j})\leq H(w)$. Moreover, since $\bar{\theta}(s)\leq\tilde{\theta}(s)$ for all $s\geq 0$, the previous inequality gives the inequality $\bar{\theta}(1/d_{i,j})\leq H(w)$. Consequently, we get $1/d_{i,j}\leq\gamma(H(w))$. Notice that since $s^{**}$ is greater than or equal to $1/\lambda$, it follows that $1/\lambda\leq\gamma(H(w))<1/L$. Hence, we obtain $d_{i,j}\geq\dfrac{1}{\gamma(H(w))}>L$. The previous inequality shows that (3.7) holds with $\rho(s):=\dfrac{1}{\gamma(s)}$ for $s\geq 0$. This concludes the proof. ◁

**Proof of Proposition 2:** Notice that by virtue of (3.7) for all $w\in\Omega$ it holds that

$$d_{i,j}\geq\rho(H(w)),\ i,j=1,...,n,\ j\neq i \tag{5.9}$$

Definition (2.4) implies that

$$\left|\dfrac{x_i-x_j}{d_{i,j}}\right|\leq 1 \text{ and } \left|\dfrac{y_i-y_j}{d_{i,j}}\right|\leq\dfrac{1}{\sqrt{p}},\text{ for all }w\in\Omega\text{ and }i,j=1,...,n,\ j\neq i \tag{5.10}$$

Moreover, since (3.2) holds and since $m\geq 2$ is the maximum number of points that can be placed within the area bounded by two concentric ellipses with semi-major axes $L$ and $\lambda$, with $L<\lambda$, and semi-minor axes $\dfrac{L}{\sqrt{p}}$ and $\dfrac{\lambda}{\sqrt{p}}$, respectively, so that each point has distance (in the metric given by (2.4)) at least $L$ from every other point, it follows that the sums $\sum_{j\neq i}V'(d_{i,j})\dfrac{(x_i-x_j)}{d_{i,j}}$, $\sum_{j\neq i}V'(d_{i,j})\dfrac{(y_i-y_j)}{d_{i,j}}$ contain at most $m$ non-zero terms, namely the terms with $d_{i,j}\leq\lambda$. Definition (3.13) in conjunction with (5.9) and (5.10) implies the following estimate for all $w\in\Omega$ and $i=1,...,n$:

$$\max\left(\left|\sum_{j\neq i}V'(d_{i,j})\frac{(x_i-x_j)}{d_{i,j}}\right|, \sqrt{p}\left|\sum_{j\neq i}V'(d_{i,j})\frac{(y_i-y_j)}{d_{i,j}}\right|\right)$$
$$\leq \max\left(\sum_{j\neq i}|V'(d_{i,j})|\left|\frac{x_i-x_j}{d_{i,j}}\right|, \sqrt{p}\sum_{j\neq i}|V'(d_{i,j})|\left|\frac{y_i-y_j}{d_{i,j}}\right|\right) \quad (5.11)$$
$$\leq \max\left(\sum_{j\neq i}b_1(\rho(H(w)))\left|\frac{x_i-x_j}{d_{i,j}}\right|, \sqrt{p}\sum_{j\neq i}b_1(\rho(H(w)))\left|\frac{y_i-y_j}{d_{i,j}}\right|\right)$$
$$\leq \sum_{j\neq i}b_1(\rho(H(w))) \leq m b_1(\rho(H(w)))$$

Definition (3.5) implies the following estimate holds for all $w \in \Omega$ and $i = 1,...,n$:

$$\frac{1}{v_{\max}\cos(\theta_i) - v^*} \leq \frac{A + (1-\cos(\varphi))H(w)}{A(v_{\max} - v^*) + (v_{\max}\cos(\varphi) - v^*)(1-\cos(\varphi))H(w)} \quad (5.12)$$

Inequality (3.16) with $R:[0,+\infty) \to \Re_+$ defined by (3.19) is a direct consequence of (5.11), (5.12) and definition (3.10).

Notice that (3.15) is equivalent to the following inequality:

$$k_i(w) \geq \max\left(\frac{1}{v^*}\sum_{j\neq i}V'(d_{i,j})\frac{(x_i-x_j)}{d_{i,j}}, -\frac{1}{v_{\max}\cos(\theta_i)-v^*}\sum_{j\neq i}V'(d_{i,j})\frac{(x_i-x_j)}{d_{i,j}}\right)$$
$$= \frac{1}{v^*}\sum_{j\neq i}V'(d_{i,j})\frac{(x_i-x_j)}{d_{i,j}} + \frac{v_{\max}\cos(\theta_i)}{v^*(v_{\max}\cos(\theta_i)-v^*)}\max\left(0, -\sum_{j\neq i}V'(d_{i,j})\frac{(x_i-x_j)}{d_{i,j}}\right) \quad (5.13)$$

Since (3.11) holds, we get for all $w \in \Omega$ and $i = 1,...,n$:

$$\max\left(0, -\sum_{j\neq i}V'(d_{i,j})\frac{(x_i-x_j)}{d_{i,j}}\right) \leq f\left(-\sum_{j\neq i}V'(d_{i,j})\frac{(x_i-x_j)}{d_{i,j}}\right). \quad (5.14)$$

Inequality (3.15) is a direct consequence of inequalities (5.13), (5.14) and definition (3.10).

Using the left inequality (3.15), the facts that $v_i \in (0, v_{\max})$, $k_i(w) > 0$, for $i = 1,...,n$, $w \in \Omega$ and definition (3.9), we get for all $w \in \Omega$ and $i = 1,...,n$:

$$F_i = -\frac{k_i(w)}{\cos(\theta_i)}(v_i\cos(\theta_i) - v^*) - \frac{1}{\cos(\theta_i)}\sum_{j\neq i}V'(d_{i,j})\frac{(x_i-x_j)}{d_{i,j}}$$
$$\geq -\frac{k_i(w)}{\cos(\theta_i)}(v_i\cos(\theta_i) - v^*) - \frac{1}{\cos(\theta_i)}k_i(w)v^* = -k_i(w)v_i \geq -k_i(w)v_{\max} \quad (5.15)$$

Using the right inequality (3.15), the facts that $v_i \in (0, v_{\max})$ for $i = 1,...,n$, $w \in \Omega$ and definition (3.9), we get for all $w \in \Omega$ and $i = 1,...,n$:

$$F_i = -\frac{k_i(w)}{\cos(\theta_i)}\left(v_i\cos(\theta_i)-v^*\right) - \frac{1}{\cos(\theta_i)}\sum_{j\neq i}V'(d_{i,j})\frac{(x_i-x_j)}{d_{i,j}}$$

$$\leq -\frac{k_i(w)}{\cos(\theta_i)}\left(v_i\cos(\theta_i)-v^*\right) + \frac{1}{\cos(\theta_i)}k_i(w)\left(v_{\max}\cos(\theta_i)-v^*\right) \quad (5.16)$$

$$= k_i(w)\left(v_{\max}-v_i\right) \leq k_i(w)v_{\max}$$

Inequalities (3.17) are direct consequences of (5.15) and (5.16).

Definition (3.9) and the fact that $v_i \in (0, v_{\max})$ for $i=1,...,n$, $w\in\Omega$, implies the following estimate for all $w\in\Omega$ and $i=1,...,n$:

$$|u_i| \leq \frac{1}{v^*}\left(\mu_1 v_{\max} + |U'(y_i)| + p\left|\sum_{j\neq i}V'(d_{i,j})\frac{(y_i-y_j)}{d_{i,j}}\right| + |F_i|\right). \quad (5.17)$$

Definition (3.14) and estimate (3.7) implies that $|U'(y_i)| \leq \kappa(H(w))$ for all $i=1,...,n$, $w\in\Omega$. Moreover, inequalities (3.17) imply that $|F_i| \leq k_i(w)v_{\max}$ for all $i=1,...,n$, $w\in\Omega$. The two previous inequalities in conjunction with (5.17) and (5.11) imply that (3.18) holds. This concludes the proof. ◁

The proof of Theorem is performed by using Barbălat's lemma ([13]) and its variant which uses uniform continuity of the derivative of a function. For the reader's convenience it is stated next.

**Lemma 3:** If a function $g \in C^2(\mathfrak{R}_+)$ satisfies $\lim_{t\to+\infty}(g(t))\in\mathfrak{R}$ and $\sup_{t\geq 0}(|\ddot{g}(t)|)<+\infty$, then, $\lim_{t\to+\infty}(\dot{g}(t))=0$.

**Proof of Theorem 1:** Using (2.3), (2.4) and (3.5), we conclude that the following equation holds for all $w\in\Omega$:

$$\nabla H(w)\dot{w} = \sum_{i=1}^{n}\left(v_i\cos(\theta_i)-v^*\right)\left(F_i\cos(\theta_i)-v_i\sin(\theta_i)u_i\right)$$

$$+\sum_{i=1}^{n}v_i\sin(\theta_i)\left(F_i\sin(\theta_i)+v_i\cos(\theta_i)u_i\right)+\sum_{i=1}^{n}U'(y_i)v_i\sin(\theta_i)$$

$$+\sum_{i=1}^{n}v_i\cos(\theta_i)\sum_{j\neq i}V'(d_{i,j})\frac{(x_i-x_j)}{d_{i,j}} \quad (5.18)$$

$$+p\sum_{i=1}^{n}v_i\sin(\theta_i)\sum_{j\neq i}V'(d_{i,j})\frac{(y_i-y_j)}{d_{i,j}} + A\sum_{i=1}^{n}\frac{\sin(\theta_i)u_i}{\left(\cos(\theta_i)-\cos(\varphi)\right)^2}$$

Using the fact that $\sum_{i=1}^{n}\sum_{j\neq i}V'(d_{i,j})\frac{(x_i-x_j)}{d_{i,j}}=0$, we get from (5.18) for all $w\in\Omega$:

$$\nabla H(w)\dot{w} = \sum_{i=1}^{n}\left(v_i \cos(\theta_i) - v^*\right)\left(F_i \cos(\theta_i) + \sum_{j\neq i} V'(d_{i,j})\frac{(x_i - x_j)}{d_{i,j}}\right)$$
$$+ \sum_{i=1}^{n} v_i \sin(\theta_i)\left(\sin(\theta_i)F_i + \left(v^* + \frac{A}{v_i\left(\cos(\theta_i) - \cos(\varphi)\right)^2}\right)u_i + U'(y_i) + p\sum_{j\neq i} V'(d_{i,j})\frac{(y_i - y_j)}{d_{i,j}}\right) \quad (5.19)$$

Finally, using (3.8), (3.9), (3.10) and the fact that $k_i(w) \geq \mu_2$ for all $i = 1,...,n$, $w \in \Omega$, we get for all $w \in \Omega$:

$$\nabla H(w)\dot{w} = -\sum_{i=1}^{n} k_i(w)\left(v_i \cos(\theta_i) - v^*\right)^2 - \mu_1 \sum_{i=1}^{n} v_i^2 \sin^2(\theta_i)$$
$$\leq -\mu_2 \sum_{i=1}^{n}\left(v_i \cos(\theta_i) - v^*\right)^2 - \mu_1 \sum_{i=1}^{n} v_i^2 \sin^2(\theta_i) \quad (5.20)$$

Let $w_0 \in \Omega$ and consider the unique solution $w(t)$ of (2.3), (3.8), (3.9), (3.10) with initial condition $w(0) = w_0$. Using the fact that the set $\Omega$ is open (recall definitions (2.6), (2.7)), we conclude that there exists $t_{\max} \in (0, +\infty]$ such that the solution $w(t)$ of (2.3), (3.9), (3.10) is defined on $[0, t_{\max})$ and satisfies $w(t) \in \Omega$ for all $t \in [0, t_{\max})$. Furthermore, if $t_{\max} < +\infty$ then there exists an increasing sequence of times $\{t_i \in [0, t_{\max}): i = 1, 2,...\}$ with $\lim_{i\to+\infty}(t_i) = t_{\max}$ and either $\lim_{i\to+\infty}(dist(w(t_i), \partial\Omega)) = 0$ or $\lim_{i\to+\infty}(|w(t_i)|) = +\infty$.

Since $w(t) \in \Omega$ for all $t \in [0, t_{\max})$, it follows from (2.6), (2.7) that $v_i(t) \in (0, v_{\max})$, $\theta_i(t) \in (-\varphi, \varphi)$ for all $t \in [0, t_{\max})$ and $i = 1,...,n$. Thus, (2.3) implies that $0 \leq \dot{x}_i(t) \leq v_{\max}$ for all $t \in [0, t_{\max})$ and $i = 1,...,n$. Moreover, inequality (5.20) gives that

$$H(w(t)) \leq H(w_0), \text{ for all } t \in [0, t_{\max}) \quad (5.21)$$

Consequently, we obtain from (3.7), (5.21) for all $t \in [0, t_{\max})$ and $i, j = 1,...,n$, $j \neq i$:

$$|\theta_i(t)| \leq \omega(H(w_0)) < \varphi, |y_i(t)| \leq \kappa(H(w_0)) < a,$$
$$x_i(0) \leq x_i(t) \leq x_i(0) + v_{\max}t, d_{i,j}(t) \geq \rho(H(w_0)) > L \quad (5.22)$$

Moreover, (2.3), (3.17) imply the following differential inequalities for all $t \in [0, t_{\max})$ and $i = 1,...,n$:

$$k_i(w(t))(v_{\max} - v_i(t)) \geq \dot{v}_i(t) \geq -k_i(w(t))v_i(t). \quad (5.23)$$

Inequalities (3.16) in conjunction with (5.21) implies that

$$k_i(w(t)) \leq M = R(H(w_0)), \text{ for all } t \in [0, t_{\max}) \text{ and } i = 1,...,n. \quad (5.24)$$

Differential inequalities (5.23) in conjunction with inequalities (5.24) imply that the following estimates hold for all $t \in [0, t_{\max})$ and $i = 1,...,n$:

$$v_i(0)\exp(-M\,t)+(1-\exp(-M\,t))v_{\max} \geq v_i(t) \geq v_i(0)\exp(-M\,t). \tag{5.25}$$

Suppose that $t_{\max} < +\infty$. Inequalities (5.21), (5.22), (5.25) and definitions (2.6), (2.7) imply that for every increasing sequence of times $\{t_i \in [0, t_{\max}) : i = 1, 2, ...\}$ with $\lim_{i \to +\infty}(t_i) = t_{\max}$ we cannot have $\lim_{i \to +\infty}(dist(w(t_i), \partial\Omega)) = 0$ or $\lim_{i \to +\infty}(|w(t_i)|) = +\infty$. Thus, we must have $t_{\max} = +\infty$.

Inequalities (3.22), (3.23) are direct consequences of inequalities (3.17), (3.18), (5.21) and (5.24). Using (2.3), (3.22), (3.23) we conclude that there exists a constant $\tilde{M} > 0$ such that

$$\frac{d}{dt}\left(\mu_2 \sum_{i=1}^{n}(v_i(t)\cos(\theta_i(t)) - v^*)^2 + \mu_1 \sum_{i=1}^{n} v_i^2(t)\sin^2(\theta_i(t))\right) \leq \tilde{M}, \text{ for all } t \geq 0. \tag{5.26}$$

Inequality (5.20) implies that

$$\mu_2 \int_0^{+\infty} \sum_{i=1}^{n}(v_i(t)\cos(\theta_i(t)) - v^*)^2 dt + \mu_1 \int_0^{+\infty} \sum_{i=1}^{n} v_i^2(t)\sin^2(\theta_i(t)) dt \leq H(w_0). \tag{5.27}$$

It follows from (5.26), (5.27) and Barbălat's lemma ([13]) that (3.20) holds.

Finally, we use Lemma 3 with $g(t) = v_i(t)$ and $g(t) = \theta_i(t)$ for $i = 1, ..., n$ in order to prove (3.21). Since (3.20) holds, it follows from Lemma 3 that it suffices to show that $\dot{u}_i(t)$ and $\dot{F}_i(t)$ are bounded. In what follows all derivatives have been calculated with the use of (2.3), (2.4), (3.8), (3.9), (3.10). We start by proving that $\dot{d}_{i,j}(t)$ is bounded for all $i, j = 1, ..., n$ with $i \neq j$. Indeed, we get:

$$\dot{d}_{i,j}(t) = \frac{(x_i(t) - x_j(t))(v_i(t)\cos(\theta_i(t)) - v_j(t)\cos(\theta_j(t)))}{d_{i,j}(t)} + p\frac{(y_i(t) - y_j(t))(v_i(t)\sin(\theta_i(t)) - v_j(t)\sin(\theta_j(t)))}{d_{i,j}(t)} \tag{5.28}$$

Using the Cauchy-Schwarz inequality and definition (2.4), we get from (5.28) for $t \geq 0$:

$$\left|\dot{d}_{i,j}(t)\right| \leq \sqrt{(v_i(t)\cos(\theta_i(t)) - v_j(t)\cos(\theta_j(t)))^2 + p(v_i(t)\sin(\theta_i(t)) - v_j(t)\sin(\theta_j(t)))^2}$$
$$\leq \sqrt{1+p}\sqrt{v_i^2(t) + v_j^2(t) - 2v_i(t)v_j(t)\cos(\theta_i(t) - \theta_j(t))} \tag{5.29}$$

Using the fact that $v_i(t) \in (0, v_{\max})$ for $i = 1, ..., n$, we get for $t \geq 0$:

$$-2v_i(t)v_j(t)\cos(\theta_i(t) - \theta_j(t)) \leq 2v_i(t)v_j(t) \leq v_i^2(t) + v_j^2(t). \tag{5.30}$$

Combining (5.29), (5.30) and using the fact that $v_i(t) \in (0, v_{\max})$ for $i = 1, ..., n$, we get for $t \geq 0$:

$$\left|\dot{d}_{i,j}(t)\right| \leq \sqrt{2(1+p)}\sqrt{v_i^2(t) + v_j^2(t)} \leq \sqrt{2(1+p)}(v_i(t) + v_j(t)) \leq 2\sqrt{2(1+p)}\,v_{\max}. \tag{5.31}$$

We next prove that $\frac{d}{dt}\left(\sum_{j\neq i}V'(d_{i,j}(t))\frac{(x_i(t)-x_j(t))}{d_{i,j}(t)}\right)$ is bounded for all $i=1,...,n$. Since (3.2) and (5.22) hold, it follows that $V'(d_{i,j}(t)), V''(d_{i,j}(t))$ are bounded for all $i,j=1,...,n$ with $i\neq j$. Since $\dot{d}_{i,j}(t)$ is bounded for all $i,j=1,...,n$ with $i\neq j$, it follows from (5.10), the fact that $d_{i,j}(t) > L$ for all $i,j=1,...,n$ with $i\neq j$, the fact that $v_i(t) \in (0, v_{max})$ for $i=1,...,n$ and the formula

$$\frac{d}{dt}\left(\sum_{j\neq i}V'(d_{i,j}(t))\frac{(x_i(t)-x_j(t))}{d_{i,j}(t)}\right)$$
$$=\sum_{j\neq i}V''(d_{i,j}(t))\dot{d}_{i,j}(t)\frac{(x_i(t)-x_j(t))}{d_{i,j}(t)}$$
$$-\sum_{j\neq i}V'(d_{i,j}(t))\dot{d}_{i,j}(t)\frac{(x_i(t)-x_j(t))}{d_{i,j}^2(t)}$$
$$+\sum_{j\neq i}V'(d_{i,j}(t))\frac{(v_i(t)\cos(\theta_i(t))-v_j(t)\cos(\theta_j(t)))}{d_{i,j}(t)}$$

that $\frac{d}{dt}\left(\sum_{j\neq i}V'(d_{i,j}(t))\frac{(x_i(t)-x_j(t))}{d_{i,j}(t)}\right)$ is bounded. Similarly, we prove that $\frac{d}{dt}\left(\sum_{j\neq i}V'(d_{i,j}(t))\frac{(y_i(t)-y_j(t))}{d_{i,j}(t)}\right)$ is bounded for all $i=1,...,n$.

We next prove that $\frac{d}{dt}(k_i(w(t)))$ is bounded for all $i=1,...,n$. Since (3.2), (5.22) and (5.10) hold, it follows that $\sum_{j\neq i}V'(d_{i,j}(t))\frac{(x_i(t)-x_j(t))}{d_{i,j}(t)}$ is bounded for all $i=1,...,n$. Using (5.21), (5.12), the fact that $\sum_{j\neq i}V'(d_{i,j}(t))\frac{(x_i(t)-x_j(t))}{d_{i,j}(t)}$, $\frac{d}{dt}\left(\sum_{j\neq i}V'(d_{i,j}(t))\frac{(x_i(t)-x_j(t))}{d_{i,j}(t)}\right)$, $u_i(t)$ are bounded and the formula

$$\frac{d}{dt}k_i(w(t))$$
$$=\frac{1}{v^*}\left(1-f'\left(-\sum_{j\neq i}V'(d_{i,j}(t))\frac{(x_i(t)-x_j(t))}{d_{i,j}(t)}\right)\right)\frac{d}{dt}\left(\sum_{j\neq i}V'(d_{i,j}(t))\frac{(x_i(t)-x_j(t))}{d_{i,j}(t)}\right)$$
$$-\frac{1}{v_{max}\cos(\theta_i(t))-v^*}f'\left(-\sum_{j\neq i}V'(d_{i,j}(t))\frac{(x_i(t)-x_j(t))}{d_{i,j}(t)}\right)\frac{d}{dt}\left(\sum_{j\neq i}V'(d_{i,j}(t))\frac{(x_i(t)-x_j(t))}{d_{i,j}(t)}\right)$$
$$+\frac{v_{max}\sin(\theta_i(t))u_i(t)}{\left(v_{max}\cos(\theta_i(t))-v^*\right)^2}f\left(-\sum_{j\neq i}V'(d_{i,j}(t))\frac{(x_i(t)-x_j(t))}{d_{i,j}(t)}\right)$$

we conclude that $\frac{d}{dt}(k_i(w(t)))$ is bounded for all $i=1,...,n$.

Definitions (2.6), (2.7) imply that $\frac{1}{\cos(\theta_i)} \leq \frac{1}{\cos(\varphi)}$ for $i=1,...,n$. Thus $\frac{1}{\cos(\theta_i(t))}$ is bounded for all $i=1,...,n$. Using the fact that $\frac{1}{\cos(\theta_i(t))}$, $\frac{d}{dt}(k_i(w(t)))$, $v_i(t)$, $k_i(w(t))$, $u_i(t)$, $F_i(t)$, $\sum_{j \neq i} V'(d_{i,j}(t)) \frac{(x_i(t) - x_j(t))}{d_{i,j}(t)}$, $\frac{d}{dt}\left(\sum_{j \neq i} V'(d_{i,j}(t)) \frac{(x_i(t) - x_j(t))}{d_{i,j}(t)}\right)$ are bounded and the formula

$$\dot{F}_i(t) = -\frac{1}{\cos(\theta_i(t))} \left(v_i(t)\cos(\theta_i(t)) - v^*\right) \frac{d}{dt}(k_i(w(t)))$$

$$-\frac{k_i(w(t))\sin(\theta_i(t))}{\cos^2(\theta_i(t))} u_i(t) \left(v_i(t)\cos(\theta_i(t)) - v^*\right)$$

$$-k_i(w(t))F_i(t) + \frac{k_i(w(t))}{\cos(\theta_i(t))} v_i(t)\sin(\theta_i(t))u_i(t)$$

$$-\frac{\sin(\theta_i(t))}{\cos^2(\theta_i(t))} u_i(t) \sum_{j \neq i} V'(d_{i,j}(t)) \frac{(x_i(t) - x_j(t))}{d_{i,j}(t)}$$

$$-\frac{1}{\cos(\theta_i(t))} \frac{d}{dt}\left(\sum_{j \neq i} V'(d_{i,j}(t)) \frac{(x_i(t) - x_j(t))}{d_{i,j}(t)}\right)$$

we conclude that $\dot{F}_i(t)$ is bounded for all $i=1,...,n$.

Combining the fact that $v_i(t) \in (0, v_{max})$ for $i=1,...,n$ with (3.20) gives the existence of a positive lower bound for all velocities, i.e., $\frac{1}{v_i(t)}$ is bounded for all $i=1,...,n$. Since (5.22) holds, it follows that $U'(y_i(t)), U''(y_i(t))$ are bounded for all $i=1,...,n$. Using the fact that $\frac{1}{v_i(t)}$, $\frac{1}{v_{max}\cos(\theta_i(t)) - v^*}$, $v_i(t)$, $u_i(t)$, $F_i(t)$, $\dot{F}_i(t)$, $U'(y_i(t))$, $U''(y_i(t))$, $\frac{d}{dt}\left(\sum_{j \neq i} V'(d_{i,j}(t)) \frac{(y_i(t) - y_j(t))}{d_{i,j}(t)}\right)$ are bounded and the formula

$$\left(v^* + \frac{A}{v_i(t)(\cos(\theta_i(t)) - \cos(\varphi))^2}\right) \dot{u}_i(t)$$

$$= \frac{F_i(t)(\cos(\theta_i(t)) - \cos(\varphi)) - 2v_i(t)u_i(t)\sin(\theta_i(t))}{v_i^2(t)(\cos(\theta_i(t)) - \cos(\varphi))^3} Au_i(t)$$

$$-p\frac{d}{dt}\left(\sum_{j \neq i} V'(d_{i,j}(t)) \frac{(y_i(t) - y_j(t))}{d_{i,j}(t)}\right) - U''(y_i(t))v_i(t)\sin(\theta_i(t))$$

$$-\left(\sin(\theta_i(t))\dot{F}_i(t) + u_i(t)\cos(\theta_i(t))F_i(t)\right)$$

$$-\mu_1\left(F_i(t)\sin(\theta_i(t)) + u_i(t)v_i(t)\cos(\theta_i(t))\right)$$

we conclude that $\dot{u}_i(t)$ is bounded for all $i=1,...,n$. This completes the proof. ◁

**Proof of Lemma 1:** Assume that $|\theta_i| \leq \varphi < \frac{\pi}{2}$ for $i = 1, 2$ with $\theta_1 \neq \theta_2$. Then, the function $f(\eta_1, \eta_2)$ in (3.26) is convex with Hessian

$$\nabla^2 f = 2 \begin{bmatrix} \cos^2(\theta_1) + p\sin^2(\theta_1) & -\cos(\theta_1)\cos(\theta_2) - p\sin(\theta_1)\sin(\theta_2) \\ -\cos(\theta_1)\cos(\theta_2) - p\sin(\theta_1)\sin(\theta_2) & \cos^2(\theta_2) + p\sin^2(\theta_2) \end{bmatrix}$$

which is a positive definite matrix since $\det(\nabla^2 f) = 4\sin^2(\theta_1 - \theta_2) > 0$. Thus, by virtue of convexity, for any $\eta_1, \eta_2 \in [0,1]$ we have

$$f(\eta_1, \eta_2) = f\big(\eta_1(1, \eta_2) + (1 - \eta_1)(0, \eta_2)\big) \leq \eta_1 f(1, \eta_2) + (1 - \eta_1) f(0, \eta_2)$$
$$f(1, \eta_2) = f\big(\eta_2(1, 1) + (1 - \eta_2)(1, 0)\big) \leq \eta_2 f(1, 1) + (1 - \eta_2) f(1, 0)$$
$$f(0, \eta_2) = f\big(\eta_2(0, 1) + (1 - \eta_2)(0, 0)\big) \leq \eta_2 f(0, 1) + (1 - \eta_2) f(0, 0)$$
$$f(\eta_1, \eta_2) \leq \eta_1 \eta_2 f(1, 1) + \eta_1(1 - \eta_2) f(1, 0) + (1 - \eta_1)\eta_2 f(0, 1) + (1 - \eta_1)(1 - \eta_2) f(0, 0)$$

which imply that
$$f(\eta_1, \eta_2) \leq \max\big(f(1,1), f(0,1), f(1,0)\big) \text{ for all } \eta_1, \eta_2 \in [0,1]. \tag{5.32}$$

Notice also that if $\theta_1 = \theta_2 = \theta \in [-\varphi, \varphi]$ it holds that

$$f(\eta_1, \eta_2) \leq \cos^2(\theta) + p\sin^2(\theta). \tag{5.33}$$

Hence, from (5.32) and (5.33) we obtain the following estimate

$$f(\eta_1, \eta_2) \leq \max(G, M) \tag{5.34}$$

where
$$G = \max\Big\{ (\cos(\theta_1) - \cos(\theta_2))^2 + p(\sin(\theta_1) - \sin(\theta_2))^2 : |\theta_i| \leq \varphi, i = 1, 2 \Big\}$$
$$M = \max\Big\{ \cos^2(\theta) + p\sin^2(\theta) : |\theta| \leq \varphi \Big\} \tag{5.35}$$

Next, define the function $g(\theta_1, \theta_2) = (\cos(\theta_1) - \cos(\theta_2))^2 + p(\sin(\theta_1) - \sin(\theta_2))^2 = g(\theta_2, \theta_1)$. Notice that

$$\frac{\partial g}{\partial \theta_1}(\theta_1, \theta_2) = -2(\cos(\theta_1) - \cos(\theta_2))\sin(\theta_1) + 2p(\sin(\theta_1) - \sin(\theta_2))\cos(\theta_1)$$
$$= 2(p-1)\cos(\theta_1)(\sin(\theta_1) - \sin(\theta_2)) + 2\sin(\theta_1 - \theta_2)$$

and consequently, $\frac{\partial g}{\partial \theta_1}(\theta_1, \theta_2) > 0$ if $-\varphi \leq \theta_2 < \theta_1 \leq \varphi$ and $\frac{\partial g}{\partial \theta_1}(\theta_1, \theta_2) < 0$ if $-\varphi \leq \theta_1 < \theta_2 \leq \varphi$. Thus, we obtain:

$$g(\theta_1, \theta_2) \leq \max\big(g(\varphi, \theta_2), g(-\varphi, \theta_2)\big) \text{ for all } (\theta_1, \theta_2) \in [-\varphi, \varphi]^2.$$

Therefore, we have:

$$g(\varphi, \theta_2) = g(\theta_2, \varphi) \leq \max\left(g(\varphi,\varphi), g(-\varphi,\varphi)\right) = g(-\varphi,\varphi) = 4p\sin^2(\varphi)$$

$$g(-\varphi, \theta_2) = g(\theta_2, -\varphi) \leq \max\left(g(\varphi,-\varphi), g(-\varphi,-\varphi)\right) = g(-\varphi,\varphi) = 4p\sin^2(\varphi)$$

Consequently, we have from the previous inequalities and (5.35) that $G = 4p\sin^2(\varphi)$. Similarly, we can prove that $M = \cos^2(\varphi) + p\sin^2(\varphi) = 1 + (p-1)\sin^2(\varphi)$. Finally, it follows from the previous calculations and (5.34) that (3.27) holds for $\eta_1, \eta_2 \in [0,1]$, $|\theta_i| \leq \varphi < \frac{\pi}{2}$ for $i=1,2$. This completes the proof. ◁

## 6. Concluding Remarks

The present work proposed decentralized control strategies for the two-dimensional movement of autonomous vehicles described by the bicycle kinematic model on lane-free roads. By leveraging appropriate tools, such as potential functions, Lyapunov functions, and barrier functions, we developed decentralized controllers that ensure that: the vehicles do not collide with each other or with the boundary of the road; the speeds of all vehicles are always positive and remain below a given speed limit; all vehicle speeds converge to a given longitudinal speed set-point; and, finally, the accelerations, lateral speeds, and orientations of all vehicles tend to zero. Future work will address the effects of nudging and appropriate notions of string-stability for vehicles operating on lane-free roads. We will also study the effect of different potential functions and the possible use of non-monotone potential functions. Finally, we will study the case where the sizes and other characteristics of the vehicles are not the same.

## Acknowledgments

The research leading to these results has received funding from the European Research Council under the European Union's Horizon 2020 Research and Innovation programme/ ERC Grant Agreement n. [833915], project TrafficFluid.

## References

[1] G. Asaithambi, V. Kanagaraj, and T. Toledo, "Driving Behaviors: Models and Challenges for Non-Lane Based Mixed Traffic", *Transportation in Developing Economies*, 2, 2016.
[2] U. Borrmann, L. Wang, A. D. Ames, and M. Egerstedt, "Control Barrier Certificates for Safe Swarm Behavior", *Proceedings of IFAC Conference on Analysis and Design of Hybrid Systems*, 2015, 68–73.
[3] K. Chavoshi, A. Kouvelas, "Cooperative Distributed Control for Lane-less and Direction-less Movement of Autonomous Vehicles on Highway Networks", *2020 - 9th Symposium of the European Association for Research in Transportation*, 2021.
[4] C. Diakaki, M. Papageorgiou, I. Papamichail, I. Nikolos, "Overview and Analysis of Vehicle Automation and Communication Systems from a Motorway Traffic Management Perspective", *Transportation Research Part A: Policy and Practice,* 75, 2015, 147-165.
[5] D. V. Dimarogonas and K. J. Kyriakopoulos, "Inverse Agreement Protocols With Application to Distributed Multi-Agent Dispersion," *IEEE Transactions on Automatic Control*, 54, 2009, 657-663.

[6] R. Delpiano, J. C. Herrera, J. Laval, J. E. Coeymans, "A Two-Dimensional Car-Following Model for Two-Dimensional Traffic Flow Problems", *Transportation Research Part C: Emerging Technologies,* 114, 2020, 504-516.

[7] T. Han and S. S. Ge, "Styled-Velocity Flocking of Autonomous Vehicles: A Systematic Design", *IEEE Transactions on Automatic Control*, 60, 2015, 2015-2030.

[8] L. Iftekhar and R. Olfati-Saber, "Autonomous Driving for Vehicular Networks with Nonlinear Dynamics," *2012 IEEE Intelligent Vehicles Symposium*, 2012, 723-729.

[9] P.A. Ioannou, and C.C. Chien, "Autonomous Intelligent Cruise Control", *IEEE Transactions on Vehicular Technology*, 42, 1993, 657-672.

[10] V. Kanagaraj, M. Treiber, "Self-driven Particle Model for Mixed Traffic and Other Disordered Flows", *Physica A: Statistical Mechanics and its Applications,* 509, 2018, 1-11.

[11] I. Karafyllis, D. Theodosis, and M. Papageorgiou, "Analysis and Control of a Non-Local PDE flow Model", *International Journal of Control*, 2020, 1–19, doi:10.1080/00207179.2020.1808902.

[12] I. Karafyllis, D. Theodosis, and M. Papageorgiou, "Nonlinear Adaptive Cruise Control of Vehicular Platoons", 2020, arXiv:2007.07054 [eess.SY], *submitted to IEEE Transactions on Automatic Control*.

[13] H. K. Khalil, *Nonlinear systems*, Prentice Hall, 2002.

[14] Y. Liang, and H. Lee, "Decentralized Formation Control and Obstacle Avoidance for Multiple Robots with Nonholonomic Constraints*", In Proceedings of American Control Conference*, 2006, 5596–5601.

[15] Y. Liu, and B. Xu, "Improved Protocols and Stability Analysis for Multivehicle Cooperative Autonomous Systems", *IEEE Transactions on Intelligent Transportation Systems*, 16, 2015 2700-2710.

[16] M. Malekzadeh, I. Papamichail, M. Papageorgiou and K. Bogenberger, "Optimal Internal Boundary Control of Lane-Free Automated Vehicle Traffic", *Transportation Research Part C,* 126, 2021, 103060.

[17] A. K. Mulla, A. Joshi, R. Chavan, D. Chakraborty, and D. Manjunath, "A Microscopic Model for Lane-less Traffic**,**" *IEEE Transactions on Control of Network Systems,* 6, 2019, 415–428.

[18] R. Olfati-Saber, "Flocking for Multi-Agent Dynamic Systems: Algorithms and Theory", *IEEE Transactions on Automatic Control*, 51, 2006, 401-420.

[19] M. Papageorgiou, K. -S. Mountakis, I. Karafyllis, I. Papamichail and Y. Wang, "Lane-Free Artificial-Fluid Concept for Vehicular Traffic," *Proceedings of the IEEE*, 109, 2021, 114-121.

[20] P. Polack, F. Altché, B. d'Andréa-Novel and A. de La Fortelle, "The Kinematic Bicycle Model: A Consistent Model for Planning Feasible Trajectories for Autonomous Vehicles?," *2017 IEEE Intelligent Vehicles Symposium (IV)*, 2017, 812-818.

[21] R. Rajamani*, Vehicle Dynamics and Control*. New York, NY, USA:Springer-Verlag, 2012.

[22] M. Rahman, M. Chowdhury, Y. Xie, and Y. He, "Review of Microscopic Lane-changing Models and Future Research Opportunities", *IEEE Transactions on Intelligent Transportation Systems*, 14, 2013, 1942–1956.

[23] H. Tanner, A. Jadbabaie, and G. J. Pappas, "Coordination of Multiple Autonomous Vehicles", *Proceedings of IEEE Mediterranean Conference on Control and Automation*, Rhodes, Greece, 2003.

[24] P. Tientrakool, Y.-C. Ho, N. Maxemchuk, "Highway Capacity Benefits From Using Vehicle-to-Vehicle Communication and Sensors for Collision Avoidance", *Proceedings of the IEEE Vehicular Technology Conference (VTC Fall)*, San Francisco, CA, 2011, 1–5.

[25] Z. Wang, G. Wu and M. J. Barth, "A Review on Cooperative Adaptive Cruise Control (CACC) Systems: Architectures, Controls, and Applications," *2018 21st International Conference on Intelligent Transportation Systems (ITSC)*, 2018, 2884-2891.

[26] P. Wieland and F. Allgöwer, "Constructive Safety Using Control Barrier Functions", *Proceedings of the 7th IFAC Symposium on Nonlinear Control Systems*, 2007, 462–467.

[27] M. T. Wolf and J. W. Burdick, "Artificial Potential Functions for Highway Driving With Collision Avoidance," *2008 IEEE International Conference on Robotics and Automation*, 2008, 3731-3736.

[28] X. Xu, J. W. Grizzle, P. Tabuada and A. D. Ames, "Correctness Guarantees for the Composition of Lane Keeping and Adaptive Cruise Control", *IEEE Transactions on Automation Science and Engineering*, 15, 2018, 1216-1229.

[29] M. M. Zavlanos, H. G. Tanner, A. Jadbabaie and G. J. Pappas, "Hybrid Control for Connectivity Preserving Flocking", *IEEE Transactions on Automatic Control*, 54, 2008, 2869-2875.

[30] Z. Zheng, "Recent Developments and Research Needs in Modeling Lane Changing", *Transportation. Research Part B: Methodological*, 60, 2014, 16–32.